\documentclass[11pt,reqno]{amsart}
\usepackage{preprint}
\usepackage{graphicx}
\usepackage{tikz}
\usetikzlibrary{arrows.meta,positioning}
\hypersetup{pdfauthor={Yuan Jiang, Zhongpu Qiu, Hongxia Hao and Jiazhen Liu},
  pdftitle={Monotonicity of the entropy production for the spatially homogeneous Landau equation with Maxwellian molecules in every dimension}}

\newcommand{\R}{\mathbb R}
\newcommand{\E}{\mathbb E}
\newcommand{\Id}{\mathrm{Id}}
\newcommand{\dd}{\mathrm d}
\newcommand{\tr}{\operatorname{tr}}
\newcommand{\diag}{\operatorname{diag}}
\newcommand{\cF}{\mathscr F}
\newcommand{\cK}{\mathcal K}
\newcommand{\cA}{\mathcal A}
\newcommand{\Eord}{\mathcal E_{\mathrm{ord}}}

\begin{document}

\title[Monotone entropy production for Landau--Maxwell]{Monotonicity of the entropy production for the\\ spatially homogeneous Landau equation\\ with Maxwellian molecules in every dimension}

\author[Y. Jiang]{Yuan Jiang}
\address[Yuan Jiang]{School of Systems Science, Beijing Normal University, Beijing 100875, China}
\author[Z. Qiu]{Zhongpu Qiu}
\address[Zhongpu Qiu]{Beijing Normal University, Beijing, China; Shanghai Artificial Intelligence Laboratory, Shanghai, China}
\email[Zhongpu Qiu]{\href{mailto:slszqzp@outlook.com}{slszqzp@outlook.com}}
\author[H. Hao]{Hongxia Hao}
\address[Hongxia Hao]{Shanghai AI Laboratory, Shanghai 200433, China}
\email[Hongxia Hao]{\href{mailto:haohongxia@pjlab.org.cn}{haohongxia@pjlab.org.cn}}
\author[J. Liu]{Jiazhen Liu}
\address[Jiazhen Liu, corresponding author]{Research Institute of Intelligent Complex Systems, Fudan University, Shanghai 200433, China}
\email[Jiazhen Liu]{\href{mailto:liujiazhen@fudan.edu.cn}{liujiazhen@fudan.edu.cn}}
\date{}


\begin{abstract}
Boltzmann's $H$-theorem asserts that the entropy of a dilute gas never decreases; the rate at which it is created, the entropy production, measures how far the gas is from thermal equilibrium. In 1966 McKean asked whether this rate itself always decreases, that is, whether a gas relaxes ever more slowly as it approaches equilibrium. Recent counterexamples show that the answer is negative for the Boltzmann equation, for the Landau equation with Coulomb interactions and for Kac's model. For the spatially homogeneous Landau equation with Maxwellian molecules, Tabary proved that the entropy production is non-increasing once the directional temperatures of the gas are close enough to each other, which always happens after an explicit waiting time, and he conjectured that no waiting time is needed.

We prove this conjecture in every dimension $d\ge2$. For every initial density that is nonnegative and has finite mass and energy, the entropy production is non-increasing from time zero on, so that the entropy is a concave function of time. No smoothness, positivity, finite entropy or finite initial entropy production is assumed: the entropy production may be infinite at time zero, in which case it becomes finite instantly and is continuous from the right at zero. The proof reduces the dynamical statement to a single inequality for smooth densities, obtained through Gaussian regularisation and an exact energy identity. This inequality is then proved by sharp functional estimates which keep track of the relations between the moments of one and the same density; they are needed precisely in the strongly anisotropic regime that was left open.
\end{abstract}

\maketitle
\preprintcontents

\section{Introduction}\label{sec:intro}

\subsection{The Landau equation and its entropy production}
In kinetic theory the state of a spatially homogeneous gas is described by the density $f=f_t(v)\ge0$ of particles with velocity $v\in\R^d$. When collisions are dominated by grazing encounters, as in a plasma, the evolution of $f$ is governed by the Landau equation. For \emph{Maxwellian molecules}, and in Tabary's normalisation~\cite{Tabary}, it reads
\begin{equation}\label{eq:landau}
 \partial_t f=Q(f):=\nabla_v\cdot\int_{\R^d}a(v-w)(\nabla_v-\nabla_w)[f(v)f(w)]\,\dd w,
 \qquad a(z)=\frac{|z|^2\Id-z\otimes z}{q},
\end{equation}
where $d\ge2$ and $q=d-1$. The matrix $a(z)$ projects onto the directions orthogonal to the relative velocity $z$. Collisions thus exchange momentum and energy between particles while conserving the total mass, momentum and energy. Multiplying $a$ by a positive constant only rescales time, so none of the results below depend on this normalisation. The weak formulation that we use is given in \cref{def:weak}.

The equilibria of~\eqref{eq:landau} are the Maxwellian (Gaussian) distributions. Relaxation towards them is measured by Boltzmann's functional $H(f)=\int f\log f$, which is minus the physical entropy. Along smooth solutions, $H$ decreases at the rate given by the \emph{entropy production}. For a density $h$ and a measurable score $p$ we write $\delta p=p(v)-p(w)$ and set
\begin{equation}\label{eq:D}
 D(h,p)=\frac12\iint h(v)h(w)\,\delta p^{\mathsf T}a(v-w)\,\delta p\,\dd v\,\dd w\;\in[0,\infty].
\end{equation}
A \emph{classical score} of an everywhere positive $h$ is a vector field $p$, measurable up to Lebesgue-null sets, such that $\log h$ is differentiable at every $v$ with gradient $p(v)$. With $p=\nabla\log f$, the quantity~\eqref{eq:D} is the rate of decrease of $H(f_t)$~\cite{VillaniMaxwell}. Since $a$ is positive semi-definite, it is non-negative, which is the $H$-theorem for this model. Moreover it vanishes exactly at Maxwellians, so it measures quantitatively the distance to equilibrium. Lower bounds for the entropy production in terms of the relative entropy are at the heart of the entropy method for the convergence to equilibrium; see~\cite{VillaniMaxwell} and the references there.

\subsection{McKean's question}
The $H$-theorem says that $H(f_t)$ decreases. A natural next question is whether it decreases in a regular way, namely whether $t\mapsto D(f_t)$ is itself non-increasing. Equivalently, one asks whether $H(f_t)$ is a convex function of time, or whether the physical entropy is concave. A positive answer means that the gas relaxes ever more slowly as it approaches equilibrium. In particular the entropy production would never exhibit a late burst after an initial phase of slow relaxation.

The question goes back to McKean~\cite{McKean}. In his study of Kac's one-dimensional caricature of a Maxwellian gas, he asked whether $H$ is even completely monotone, that is, whether $(-1)^n\frac{\dd^n}{\dd t^n}H(f_t)\ge0$ for all $n\ge1$. Lieb~\cite{Lieb} and Olaussen~\cite{Olaussen} showed that complete monotonicity fails for the Boltzmann equation, but their arguments concern derivatives of high order. The case $n=2$, that is, the monotonicity of the entropy production itself, remained open for a long time.

The picture has changed recently. Silvestre~\cite{Silvestre} constructed solutions of the Boltzmann equation, for hard spheres and for Maxwell molecules, along which the entropy production increases for some time. Junn\'e, Winter and Yolda\c s~\cite{JWY} obtained a counterexample for the Landau equation with Coulomb interactions, and Wang~\cite{Wang} one for Kac's equation. On the positive side, it has long been known that the Fisher information $\int|\nabla f|^2/f$, a simpler relative of~\eqref{eq:D}, is non-increasing along the Landau equation with Maxwellian molecules~\cite{VillaniMaxwell}. Guillen and Silvestre~\cite{GuillenSilvestre} have extended this to a wide range of interactions, including the Coulomb case, which is the key to their global regularity theory. For radially symmetric data the Landau--Maxwell equation reduces to a linear Fokker--Planck equation~\cite[Remark~2.3]{Tabary}, and~\eqref{eq:D} becomes a Fisher information, which is again non-increasing. The difficulty is therefore entirely due to anisotropy: the entropy production couples the second derivatives of $\log f$ to the differences between the directional temperatures of the gas.

\subsection{Tabary's theorem and the anisotropic layer}
Normalise $f$ to unit mass and zero mean velocity, and write its covariance matrix in principal axes as $T=\diag(T_1,\dots,T_d)$ with $\tr T=d$. The eigenvalues $T_i$ are the \emph{directional temperatures}. For Maxwellian molecules they obey a closed linear equation, and relax exponentially to the common value $1$,
\[
 T_i(t)-1=e^{-4dt/q}\bigl(T_i(0)-1\bigr).
\]
Tabary~\cite{Tabary} recently derived an explicit formula for $\frac{\dd}{\dd t}D(f_t)$ and showed that it is non-positive as soon as $T_{\max}=\max_iT_i\le1+C_d$, for an explicit constant $C_d>0$; for instance $C_2=(13-\sqrt{133})/18\approx0.08$ and $C_3=(20-2\sqrt{70})/15\approx0.22$. Since the temperatures relax at a known rate, this gives monotonicity of $D(f_t)$ for all $t\ge t_0$, where
\[
 t_0=\frac{q}{4d}\,\log\frac{T_{\max}(0)-1}{C_d},
\]
and $t_0=0$ if the right-hand side is negative~\cite[Theorem~1.1]{Tabary}. Tabary also explained why the remaining \emph{anisotropic layer} $[0,t_0)$ lies beyond his method. His argument discards several favourable terms in the formula for $\frac{\dd}{\dd t}D$, and these terms are indispensable when one direction is much hotter than another~\cite[Remark~4.3]{Tabary}. He conjectured that monotonicity holds for all initial data from time zero~\cite[discussion after Remark~1.2]{Tabary}.

\begin{figure}[t]
 \centering
 \includegraphics[width=\textwidth]{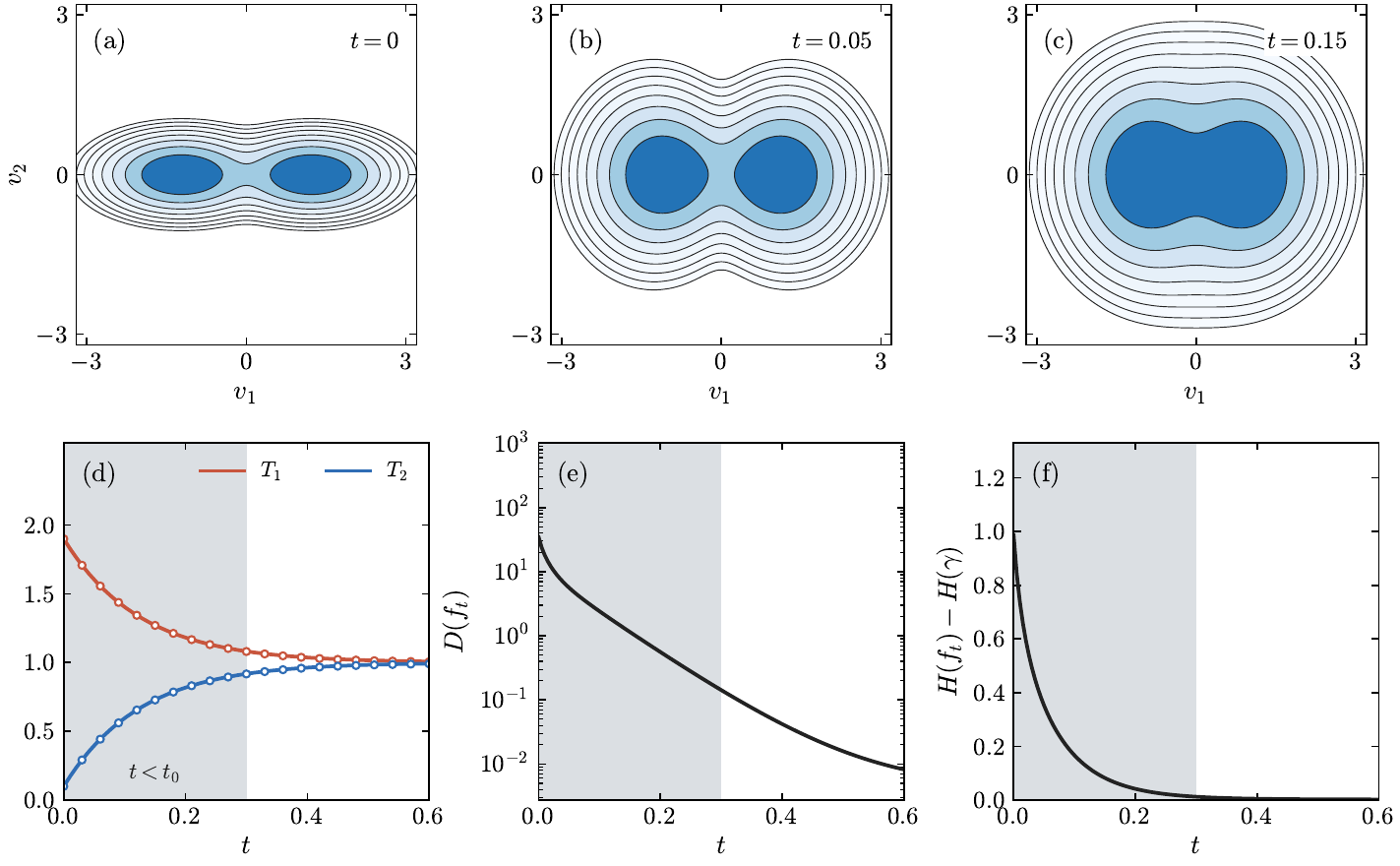}
 \caption{A numerical illustration in dimension $d=2$. The initial density is an equal mixture of two Gaussians centred at $(\pm\sqrt{3/2},0)$ with covariance $\diag(0.4,0.1)$, so that $T(0)=\diag(1.9,0.1)$. (a)--(c) The density $f_t$ at $t=0$, $0.05$ and $0.15$, with level sets at $2^{-k}\max f_t$. (d) The directional temperatures $T_1,T_2$; the markers are computed from the numerical solution and the curves are the exact law $T_i(t)-1=e^{-8t}(T_i(0)-1)$. (e) The entropy production $D(f_t)$, on a logarithmic scale. (f) The entropy $H(f_t)-H(\gamma)$ relative to the standard Gaussian $\gamma$, which is convex in $t$. In (d)--(f) the shaded strip is the anisotropic layer $[0,t_0)$, $t_0\approx0.30$, not covered by Tabary's theorem. The solution was computed by finite differences; the figure is an illustration and is not used in any proof.}
 \label{fig:illustration}
\end{figure}

\Cref{fig:illustration} shows a typical solution in dimension two whose initial temperatures are $T(0)=\diag(1.9,0.1)$. The initial density consists of two elongated bumps. Within a fraction of a time unit it becomes nearly isotropic, and the entropy production drops by two orders of magnitude. Most of this drop takes place inside the layer $[0,t_0)$, where Tabary's criterion does not apply.

\subsection{Main result}
We prove Tabary's conjecture in every dimension, and for every physically admissible initial datum. An \emph{admissible density} is nonnegative almost everywhere on $\R^d$, with $f,\ |v|^2 f\in L^1(\dd v)$. For such densities the formula~\eqref{eq:D} need not make sense, because $\log f$ need not be differentiable. We therefore use its closed extension $D_{\rm cl}$ with values in $[0,\infty]$, defined in \cref{def:closedD} through the weak derivatives of $\sqrt f$. It coincides with~\eqref{eq:D} whenever the latter is meaningful (\cref{lem:closedD}).
\begin{theorem}\label{thm:main}
Let $(f_t)_{t\ge0}$ be any weak solution with admissible initial density $f_{\rm in}$ and mass $m=\int f_{\rm in}$. Then
\begin{equation}\label{eq:main}
 D_{\rm cl}(f_t)\le D_{\rm cl}(f_s)\quad(0\le s\le t),\qquad
 \lim_{t\downarrow0}D_{\rm cl}(f_t)=D_{\rm cl}(f_{\rm in})\quad\text{in }[0,\infty].
\end{equation}
If $m>0$, for every $t>0$ there is a representative $h_t=f_t$ almost everywhere with a classical score $p_t$. Every such pair gives $D(h_t,p_t)=D_{\rm cl}(f_t)<\infty$. Thus $D(h_t,p_t)\le D(h_s,p_s)$ for every $0<s\le t$ and every choice of such pairs at the two times. If $m=0$, every slice vanishes almost everywhere and its closed production is zero.
\end{theorem}

Some comments are in order.
\begin{itemize}
\item No smoothness, positivity, finite entropy, finite Fisher information or higher moments are assumed. The initial entropy production may be infinite. It then becomes finite at every positive time, and~\eqref{eq:main} shows that it is continuous from the right at $t=0$ in $[0,\infty]$.
\item The weak solutions of \cref{def:weak} are unique (\cref{prop:evolution}). In \cref{sec:discussion} we check that the classical Cauchy solution of Villani~\cite{VillaniMaxwell}, which is the solution used by Tabary, is one of them. The positive-time part of \cref{thm:main} is therefore exactly Tabary's conjecture.
\item Since $\frac{\dd}{\dd t}H(f_t)=-D(f_t)$ along these solutions, \cref{thm:main} says that $H(f_t)$ is convex in time, and the physical entropy is concave.
\item The theorem compares values of the entropy production, not derivatives. In particular, we do not claim that $D(f_t)$ has a finite right derivative at $t=0$ for rough data.
\end{itemize}

\subsection{Strategy of the proof}
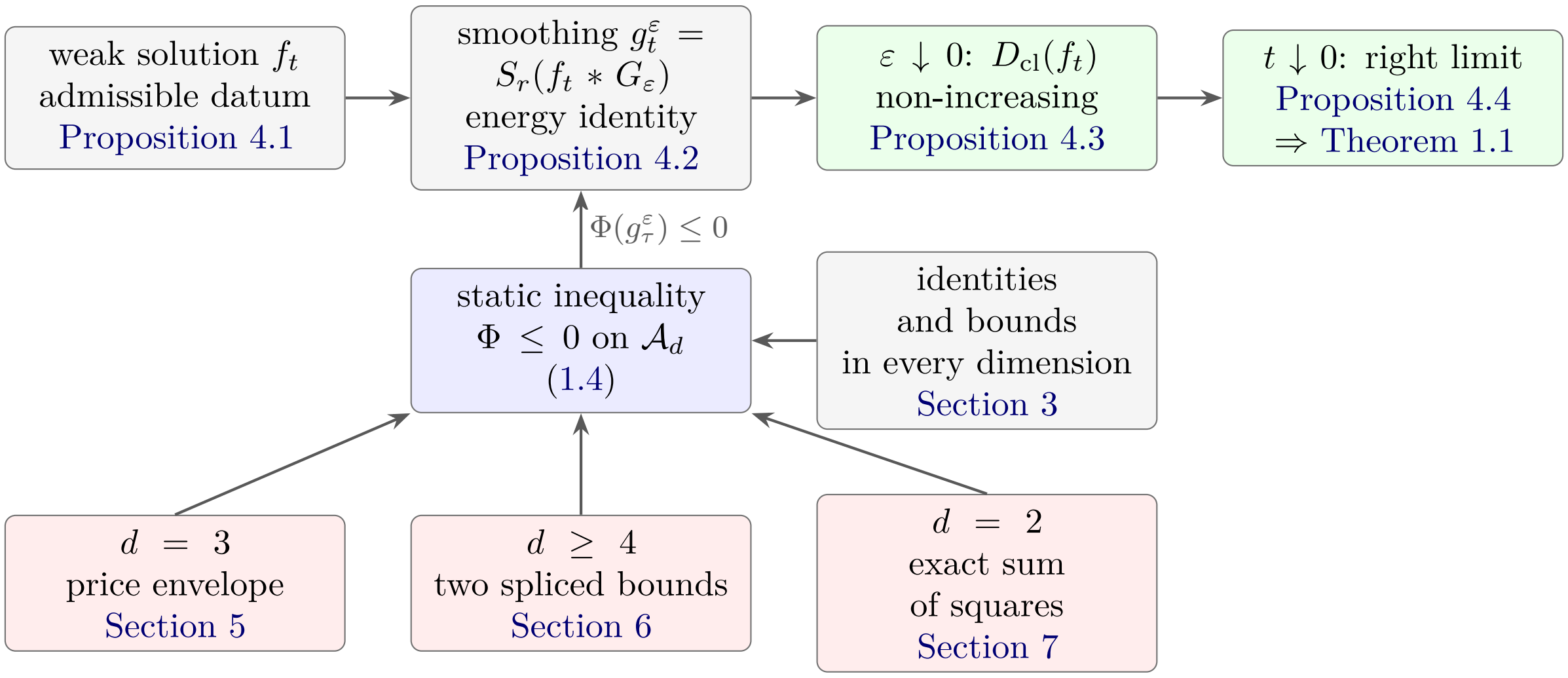
\begin{figure}[t]
\centering
\begin{tikzpicture}[
  box/.style={draw=black!55, rounded corners=3pt, align=center, font=\small, inner sep=4pt,
              minimum height=1.3cm, text width=3.15cm, fill=#1},
  box/.default=white,
  arr/.style={-{Stealth[length=2.2mm]}, thick, draw=black!65},
  lab/.style={font=\footnotesize, text=black!65, fill=white, inner sep=1pt}]
 \node[box=black!4]  (weak) at (0,0)    {weak solution $f_t$\\ admissible datum\\ \cref{prop:evolution}};
 \node[box=black!4]  (conv) at (4.1,0) {smoothing $g^\varepsilon_t=S_r(f_t*G_\varepsilon)$\\ energy identity\\ \cref{prop:gaussian-comparison}};
 \node[box=green!8]  (mono) at (8.2,0)  {$\varepsilon\downarrow0$: $D_{\rm cl}(f_t)$\\ non-increasing\\ \cref{prop:graph-comparison}};
 \node[box=green!8]  (main) at (12.3,0){$t\downarrow0$: right limit\\ \cref{prop:initial}\\ $\Rightarrow$ \cref{thm:main}};
 \node[box=blue!8]   (static) at (4.1,-2.45) {static inequality\\ $\Phi\le0$ on $\cA_d$\\ \eqref{eq:static}};
 \node[box=black!4]  (ident)  at (8.2,-2.45)  {identities and bounds\\ in every dimension\\ \cref{sec:identity}};
 \node[box=red!7]    (d3) at (0,-4.9)    {$d=3$\\ price envelope\\ \cref{sec:d3}};
 \node[box=red!7]    (d4) at (4.1,-4.9) {$d\ge4$\\ two spliced bounds\\ \cref{sec:d4}};
 \node[box=red!7]    (d2) at (8.2,-4.9)  {$d=2$\\ exact sum of squares\\ \cref{sec:d2}};
 \draw[arr] (weak) -- (conv);
 \draw[arr] (conv) -- (mono);
 \draw[arr] (mono) -- (main);
 \draw[arr] (static) -- node[lab, right=1pt] {$\Phi(g^\varepsilon_\tau)\le0$} (conv);
 \draw[arr] (ident) -- (static);
 \draw[arr] (d3.north) -- (static.south west);
 \draw[arr] (d4) -- (static);
 \draw[arr] (d2.north) -- (static.south east);
\end{tikzpicture}
\caption{Structure of the proof. The upper row is the dynamical argument: the weak solution is smoothed by a Gaussian, the Fisher energy of the smoothed solution obeys an exact identity whose integrand is $\Phi$, and the smoothing and the initial time are then removed. The static inequality $\Phi\le0$ (middle row) is the only input from the lower row, where it is proved separately in the three dimensional regimes using the common identities of \cref{sec:identity}.}
\label{fig:roadmap}
\end{figure}

The proof has two parts of a different nature, summarised in \cref{fig:roadmap}.

\emph{From the dynamics to a static inequality.} For Maxwellian molecules the coefficient $a*f$ of~\eqref{eq:landau} depends on $f$ only through its mass, mean and covariance, which evolve explicitly. Once the initial covariance is known, \eqref{eq:landau} therefore becomes a \emph{linear} Fokker--Planck equation~\eqref{eq:linear} whose coefficients are prescribed in advance. Its solutions remain solutions after convolution with a Gaussian and a dilation. For the smoothed solution $g^\varepsilon_t$ we prove the exact energy identity
\[
 \cF(g^\varepsilon_t)-\cF(g^\varepsilon_s)=\int_s^t\Phi(g^\varepsilon_\tau)\,\dd\tau,
\]
where $\cF=D_{\rm cl}+d$ is the Fisher energy of \cref{def:closedD} and $\Phi(f)$ is an explicit functional of the logarithmic derivatives of $f$ up to second order (\cref{sec:identity}). Thus \cref{thm:main} follows from the \emph{static inequality}
\begin{equation}\label{eq:static}
 \Phi(f)\le0\qquad(f\in\cA_d),
\end{equation}
where $\cA_d$ is the class of smooth positive densities of \cref{def:staticclass}. This inequality no longer involves time. The passage back to rough data uses lower semicontinuity of the Fisher energy, recovery of the energy under Gaussian smoothing, and a priori bounds that make the entropy production finite at every positive time (\cref{sec:bridge}).

\begin{figure}[t]
 \centering
 \includegraphics[width=\textwidth]{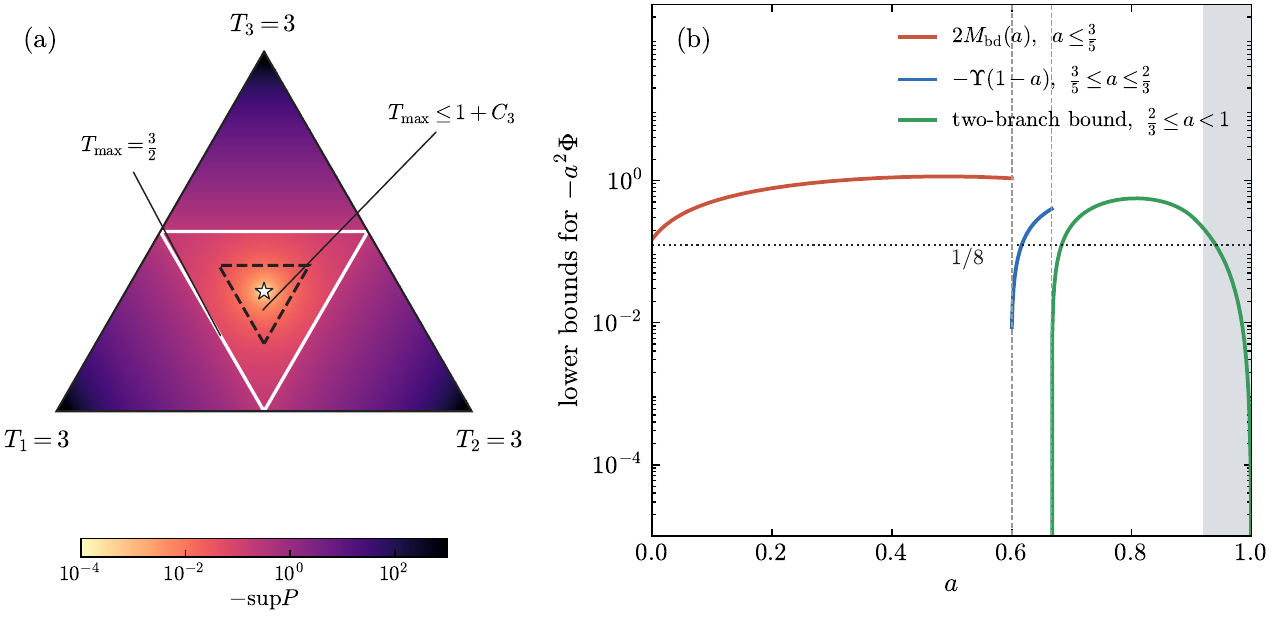}
 \caption{Certified margins in the static inequality. (a) Dimension three. Each point of the triangle is a temperature vector $(T_1,T_2,T_3)$ with $T_1+T_2+T_3=3$, and the colour is the margin $-\sup_IP(T,I)$ guaranteed by \cref{lem:d3moderate,lem:d3hot}, where $P\ge\Phi$ is the price envelope of \cref{sec:d3}. Inside the white triangle $T_{\max}\le3/2$ the margin is $\sum_i(T_i-1)^2/6$; outside it is $\frac{135}{64}(\frac{T_{\max}-1}{3-T_{\max}})^2$. The dashed triangle is Tabary's region $T_{\max}\le1+C_3$ and the star marks equilibrium. (b) Dimension two. Lower bounds for $-a^2\Phi$ as functions of the smaller temperature $a$, given by \cref{thm:d2small} for $a\le3/5$ and by the two cases in the proof of \cref{thm:d2moderate} for $3/5\le a<1$. The two-branch bound vanishes at $a=2/3$ and $a=1$; at these points the estimate $\Phi\le-J^2/2$ of \cref{thm:d2moderate} still applies. The dotted line marks $1/8$, and the shaded strip is Tabary's region $a\ge1-C_2$.}
 \label{fig:certificates}
\end{figure}

\emph{The static inequality.} The heart of the paper is the proof of~\eqref{eq:static} (\cref{sec:identity,sec:d3,sec:d4,sec:d2}). In terms of the directional temperatures $T_i$ and the directional Fisher informations $I_i=\int f(\partial_i\log f)^2$, the functional $\Phi$ consists of a positive part, linear in the $I_i$ with coefficients that grow with the anisotropy, and of a negative part given by squares of second derivatives of $\log f$ along translations and rotations. Near isotropy the negative part wins easily; this is essentially the regime covered by Tabary. In the strongly anisotropic regime, however, the natural lower bounds for the squares are not sufficient if the various moments of $f$ are treated as independent parameters. The key point of our approach is to retain enough of the relations that hold between the moments of one and the same density. These relations come from integration by parts and from the geometry of the rotation fields on spheres. How this is done depends on the dimension.
\begin{itemize}
\item In dimension $d=3$ (\cref{sec:d3}) we bound $\Phi$ by a function of finitely many scalar moments. All angular quantities are then eliminated at once by an explicit one-dimensional optimisation, which we call a price envelope. The envelope is shown to be negative by a case analysis on the temperatures. The hard case is a single hot direction, where we use a concavity argument in the cold temperatures and a monotone elimination of the coldest Fisher information.
\item In dimensions $d\ge4$ (\cref{sec:d4}), two different upper bounds are used in two complementary regions of the temperature simplex. The first one applies when all temperatures are at most $d-1$. The second one covers the extreme case in which a single direction carries almost all of the energy.
\item In dimension $d=2$ (\cref{sec:d2}) there is only one rotation field, and the margin is smallest. We write $\Phi$ exactly as an explicit function of a few moments plus a sum of non-negative squares. Two mixed moments are fixed by integration by parts in terms of the others. Keeping these relations, a single multiplier and elementary Fenchel inequalities give a uniform negative bound when the smaller temperature $a$ is at most $3/5$. A separate completion of squares covers $3/5\le a\le1$.
\end{itemize}
The resulting certified margins are displayed in \cref{fig:certificates}; they are quantitative and degenerate only at equilibrium.

\subsection*{Organisation}
\Cref{sec:prelim} collects the normalisation, the notion of weak solution and the closed entropy production. \Cref{sec:identity} derives the identities and bounds that are common to all dimensions. \Cref{sec:bridge} constructs the weak solution and reduces \cref{thm:main} to the static inequality~\eqref{eq:static}. \Cref{sec:d3,sec:d4,sec:d2} prove~\eqref{eq:static} in dimension three, in dimensions four and higher, and in dimension two, respectively. Finally, \cref{sec:discussion} compares our solutions with the classical Cauchy solution used by Tabary and discusses an open question.

\section{Preliminaries}\label{sec:prelim}

\subsection{Symmetries and normalisation}
For $m>0$ let $c=m^{-1}\int vf_{\rm in}$ and $\sigma^2=(dm)^{-1}\int |v-c|^2 f_{\rm in}$. A positive-mass Lebesgue density gives positive variance in every direction. One fixed orthogonal matrix $O$ gives a normalised weak solution $g_r(v)=\sigma^d f_{r/m}(c+\sigma Ov)/m$. For every $r\ge0$,
\begin{equation}\label{eq:normal}
 \int g_r=1,\qquad \int vg_r=0,\qquad T(r)=\int (v\otimes v)g_r=\diag(T_1(r),\ldots,T_d(r)),\qquad \tr T(r)=d.
\end{equation}
The kernel $a(z)$ is homogeneous of degree two, so velocity dilation cancels the two gradient factors in $\cF$. Translation and orthogonal changes also preserve $\cF$. The mass factor is quadratic in both the collision operator and the graph energy, giving $\cF(f_t)=m^2\cF(g_{mt})$ and $D_{\rm cl}(f_t)=m^2D_{\rm cl}(g_{mt})$, including infinite values. Testing with cutoffs of $v_i v_j$ gives, with a right derivative at $r=0$,
\begin{equation}\label{eq:Tode}
 \dot T=\frac{4d}{q}(\Id-T),\qquad T(r)=\Id+e^{-4dr/q}(T(0)-\Id).
\end{equation}
Write $f$ for the normalised density and $M=(d\Id-T)/q$. Each $M_i=(d-T_i)/q$ is bounded below by $\min\{1,M_i(0)\}>0$ throughout the evolution. Expanding $a*f=a+M$ gives the normalised forward equation
\begin{equation}\label{eq:linear}
 \partial_t f=\nabla\cdot\bigl((a+M)\nabla f+vf\bigr)
\end{equation}
in distributions.

\subsection{Score quantities}
For a positive classical slice in \eqref{eq:normal} of finite graph energy, write $u=\log f$, $L_{ij}=v_i\partial_j-v_j\partial_i$ and $\E h=\int h f$. With $I_i=\E[(\partial_i u)^2]$ and $\tau_{ij}=\E[(L_{ij}u)^2]$, the identity $q\,s^{\mathsf T}a(v)s=\sum_{i<j}(v_i s_j-v_j s_i)^2$ gives
\begin{equation}\label{eq:Dscalar}
 D(f)=\int f\,(\nabla u)^{\mathsf T}(a+M)\nabla u-d
     =\sum_i M_iI_i+q^{-1}\sum_{i<j}\tau_{ij}-d .
\end{equation}
Here $D(f)=D(f,\nabla u)$, by \cref{lem:closedD}.

\begin{definition}[Class for the static inequality]\label{def:staticclass}
A strictly positive admissible density belongs to $\cA_d$ if it has unit mass, total second moment $d$, finite graph energy $\cF(f)$ and smooth $u=\log f$. In dimension two require also $\int vf=0$. The fields $\partial_{ij}u$, $L_{ij}u$, $\partial_k L_{ij}u$ and $L_{ij}L_{k\ell}u$ belong to $L^2(f\,\dd v)$ for all indicated indices, taking $i<j$ and $k<\ell$ for rotation pairs. Here $T_i=\E[v_i^2]$ and $M_i=(d-T_i)/q$ in the chosen coordinates.
\end{definition}

\subsection{Weak solutions and the closed dissipation}
Write $\langle f,\varphi\rangle=\int f\varphi$ and $B_\varphi(f)=\iint f(v)f(w)K_\varphi(v,w)\,\dd v\,\dd w$, where
\[
 K_\varphi(v,w)=a(v-w):\nabla^2\varphi(v)-2(v-w)\cdot\nabla\varphi(v).
\]
\begin{definition}[Weak solutions]\label{def:weak}
A weak solution has $f_0=f_{\rm in}$ almost everywhere and admissible slices for every $t\ge0$. It is jointly measurable up to a product-null set on $[0,\infty)\times\R^d$ and has locally bounded total second moment there in time. Each pairing $\langle f_t,\varphi\rangle$, $\varphi\in C_c(\R^d)$, is continuous on $[0,\infty)$. For every $\varphi\in C_c^\infty(\R^d)$, $B_\varphi(f_t)$ converges absolutely at each $t\ge0$. For every such $\varphi$ and $\eta\in C_c^\infty((0,\infty))$ the following integrand belongs to $L^1([0,\infty))$ and
\[
 \int_0^\infty\bigl[\eta'(t)\langle f_t,\varphi\rangle+\eta(t)B_\varphi(f_t)\bigr]\,\dd t=0.
\]
\end{definition}

\begin{definition}[Closed dissipation]\label{def:closedD}
For an admissible $f$ of mass $m$, take a nonnegative representative, put $\varrho=\sqrt f$ and set
\begin{equation}\label{eq:closedgraph}
 \cF(f)=\inf_{g}4\iint f(w)\,g(v)^{\mathsf T}a(v-w)g(v)\,\dd v\,\dd w,
 \qquad D_{\rm cl}(f)=\cF(f)-dm^2.
\end{equation}
The infimum runs over all locally integrable distributional gradients $g$ of $\varrho$ and is $+\infty$ if none exists. Both functionals depend only on the almost-everywhere class of $f$.
\end{definition}
For a centred probability density let $T=\int(v\otimes v)f$ and $M_f=(\tr T\,\Id-T)/q$. Expanding the quadratic kernel gives
\[
 \cF(f)=4\int (\nabla\varrho)^{\mathsf T}M_f\nabla\varrho+\frac4q\sum_{i<j}\int|L_{ij}\varrho|^2
\]
when $\varrho\in H^1$ and its weak rotation derivatives belong to $L^2$. These conditions characterise finite energy. Indeed, for a unit vector $e$, $q\,e^{\mathsf T}M_fe=\int |v-(v\cdot e)e|^2f>0$ since a line has Lebesgue measure zero. Thus $M_f$ controls the full gradient. For arbitrary $m>0$, the graph energy of $f$ is $m^2$ times this expression for the centred probability density $f(c+\cdot)/m$.

\begin{lemma}\label{lem:closedD}
Every admissible $f$ satisfies $\cF(f)\ge dm^2$. If $\cF(f)<\infty$, a weak gradient $g$ of $\sqrt f$ exists. For every such $g$, put $s=2g/\sqrt f$ on $\{f>0\}$ and $s=0$ elsewhere. Then $D(f,s)=D_{\rm cl}(f)$, also the value for every classical representative and its score. At fixed mass and momentum, both $\cF$ and $D_{\rm cl}$ are lower semicontinuous under $L^1$ convergence of the densities and convergence of their total second moments.
\end{lemma}
\begin{proof}
At finite energy the Sobolev derivative $g$ vanishes almost everywhere on $\{\sqrt f=0\}$ by the level-set property of weak gradients. Hence the diagonal terms in the polarisation of \eqref{eq:D} add to $\cF(f)$. The flux $J=fs=2\sqrt f\,g=\nabla f$ satisfies $|v|J\in L^1$ by Cauchy--Schwarz. The cross term $J(v)^{\mathsf T}a(v-w)J(w)$ is absolutely integrable by Cauchy--Schwarz for $a$. Fix a smooth radial cutoff $0\le\chi\le1$ equal to one on $B_1$ and zero outside $B_2$, and put $\chi_R(v)=\chi(v/R)$. Insert $\chi_R(v)\chi_R(w)$ before expanding the kernel. Integration by parts gives $\int\chi_R J=O(m/R)$ and $\int\chi_R v_jJ_i\to-m\delta_{ij}$. For $|P(v)|\le |v|^2$, dominated convergence gives $R^{-1}\!\int\chi_R PJ\to0$, since $R^{-1}\chi_R|P J|\le2|v||J|$. Thus the products of a quadratic and a zeroth moment vanish. The products of first moments contribute $(-2d+d^2+d)m^2/q=dm^2$. Removing both cutoffs in the full cross term gives $D(f,s)=\cF(f)-dm^2\ge0$. The lower bound at infinite energy follows at once.

For a classical representative the chain rule gives the pointwise derivative of its square root. Gradients of its mollifications converge pointwise to this derivative and in $L^2$ to the weak gradient. Hence $g=\tfrac12\sqrt h\,\nabla\log h$ almost everywhere for every classical representative $h$.

For lower semicontinuity, centre and normalise at positive mass. Truncating $|v|^2$ shows that $L^1$ convergence together with convergence of total second moments gives convergence of the moment matrices. On a subsequence realising a finite lower limit, coercivity bounds the ordinary and rotation gradients in $L^2$. Strong $L^2$ convergence of the square roots identifies the weak limits as derivatives of $\sqrt f$. The limiting quadratic form is weakly lower semicontinuous. At zero mass all energies vanish.
\end{proof}

\section{Identities valid in every dimension}\label{sec:identity}

Let $d\ge2$ and $q=d-1$. Assume $f$ is a strictly positive admissible density of unit mass and total second moment $d$, with the regularity and integrability of \cref{def:staticclass}. Write $\E$ for integration against $f$, $u=\log f$, $s=\nabla u$ and $H=\nabla^2u$.

The coordinate second moments $T_i=\E[v_i^2]$ are positive and sum to $d$. Thus $M_i=(d-T_i)/q>0$ and $\sum_i M_i=d$. After centring $f$, the positive definite matrix in the expansion of \eqref{eq:closedgraph} controls $\nabla\sqrt f$ in $L^2$. Thus $fs_i^2=4|\partial_i\sqrt f|^2$ is integrable. For $X=\partial_i$ or $L_{ij}$, integration of $X(fg)$ against a radial cutoff gives $\E[Xg]=-\E[gXu]$ whenever $g,Xg,gXu\in L^1(f)$. The Cartesian boundary term is bounded by $C\E|g|/R$ and the angular one vanishes. Taking $X=\partial_i$, $g=v_i$ gives $\E[v_i s_i]=-1$, hence $I_i=\E[s_i^2]\ge1/T_i$ by Cauchy--Schwarz.
For a plane $p=(i,j)$, $i<j$, set $L_p=v_i\partial_j-v_j\partial_i$ and put
\[
 w_{ij}=L_pu,\quad \tau_{ij}=\E[w_{ij}^2],\quad C_{ij}=\E[H_{ij}w_{ij}],\quad
 \delta_{ij}=M_i-M_j,\quad \Delta_{ij}=I_i-I_j .
\]
With $J=q^{-1}\sum_{i<j}\tau_{ij}$ and $\cK=2q^{-2}\sum_{p,r}\E[(L_p L_r u)^2]$, where $p,r$ run independently over the coordinate planes, the spatial form is
\begin{equation}\label{eq:phi}
\begin{aligned}
 \Phi(f)={}&\sum_i b_i I_i-\cK-\frac4q\sum_{k,\,i<j}M_k\E[(\partial_k w_{ij})^2]
               -2\sum_{i,j}M_iM_j\E[H_{ij}^2]\\
 &-\frac4q\sum_{i<j}\delta_{ij}(\Delta_{ij}-C_{ij}),\qquad b_i=\frac{6d-(4d+2)M_i}{q}.
\end{aligned}
\end{equation}
Along a classical solution, $\Phi(f_t)$ is the time derivative of the classical Fisher integral at each $t>0$ for which differentiation under the integral and the Green formulas are justified. The corresponding integrated identity for Gaussian convolutions uses $\Phi$ at almost every time.

\begin{lemma}[Spherical estimate]\label{lem:sphere}
$\displaystyle \cK\ge\frac2q\bigl(J^2+(d-2)J\bigr)$.
\end{lemma}
\begin{proof}
Set $A=|v|^2\Id-vv^{\mathsf T}$, $c=v\cdot s$ and $\Theta=\sum_p L_p^2u$. At $v=\rho e_1$, the only nonzero rows of $(L_pL_ru)_{p,r}$ can have $p=(1,j)$. For $j,k,\ell>1$ their entries are
\[
 L_{1j}L_{1k}u=\rho^2 H_{jk}-\delta_{jk}c,\qquad
 L_{1j}L_{k\ell}u=\rho(\delta_{jk}s_\ell-\delta_{j\ell}s_k)\quad(k<\ell),
\]
where the deltas in this display are Kronecker symbols. In the second block each $\rho^2s_k^2$ occurs $d-2$ times in the squared norm. An orthogonal change of coordinates acts orthogonally on the plane generators (by its action on $\bigwedge^2\R^d$), preserving the trace and squared norm of this matrix. Thus at every $v$, $\Theta=\tr(HA)-qc$ and
\[
 \sum_{p,r}(L_pL_ru)^2-(d-2)\sum_p(L_pu)^2
     =\tr(HAHA)-2c\tr(HA)+qc^2 .
\]
For $v\ne0$, $A/|v|^2$ projects onto the $q$-dimensional space $v^\perp$. Cauchy--Schwarz there gives $(\tr(HA))^2\le q\tr(HAHA)$, also valid at $v=0$ where $A=0$. Thus $\Theta^2$ is at most $q$ times the displayed expression.

Angular integration by parts gives $\E[\Theta]=-qJ$. Integrating the pointwise bound and using $(\E[\Theta])^2\le\E[\Theta^2]$ proves the estimate.
\end{proof}

For $X=\partial_i$ or $L_p$ put $z=Xu$. Finite graph energy and the static hypotheses give $z,Xz\in L^2(f)$. With $g_\varepsilon(z)=z^3/(1+\varepsilon z^2)$, weighted integration by parts gives $\E[zg_\varepsilon(z)]=-\E[g_\varepsilon'(z)Xz]$. Indeed, for fixed $\varepsilon>0$ the test grows at most linearly and has bounded derivative. Use a smooth radial cutoff equal to one on $|v|\le R$, with gradient bounded by $C/R$ and angular derivative zero. Since $(g_\varepsilon')^2\le9zg_\varepsilon$, Cauchy--Schwarz, followed by monotone convergence as $\varepsilon\downarrow0$, gives the first estimate below. The bound $|g_\varepsilon'(z)|\le3z^2$ then gives the second identity by dominated convergence.
\begin{equation}\label{eq:cubictest}
 \E z^4\le9\E[(Xz)^2],\qquad 3\E[z^2Xz]=-\E z^4.
\end{equation}

\begin{lemma}[Reduction in one coordinate plane]\label{lem:pairbudget}
Put $B_d=\sum_i(b_i I_i-2M_i^2 I_i^2)$. For each $i<j$, abbreviate $W=w_{ij}$, $C=C_{ij}$, $\Delta=\Delta_{ij}$, $\tau=\tau_{ij}$, $\mu_i=M_i$, $\mu_j=M_j$, and set
\[
 \alpha=\frac4q\left(\frac{\mu_i}{I_j}+\frac{\mu_j}{I_i}\right),\qquad e=\frac{2(\mu_i-\mu_j)}q,\qquad A=4\mu_i\mu_j.
\]
Then $\Phi\le B_d+\sum_{i<j}Q_{ij}-\cK$, where
\[
 Q_{ij}=-\frac4q\bigl(\mu_i\E[(\partial_i W)^2]+\mu_j\E[(\partial_j W)^2]\bigr)-A\E[H_{ij}^2]-2e(\Delta-C).
\]
For every $w\in\R$,
\begin{equation}\label{eq:pairdual}
 Q_{ij}\le\frac{(w-e)^2}{\alpha}-2w\Delta+\frac{w^2\tau}{A} .
\end{equation}
\end{lemma}
\begin{proof}
Integration by parts gives $\E[H_{ii}]=-I_i$, hence $\E[H_{ii}^2]\ge I_i^2$. Each mixed Hessian occurs twice in \eqref{eq:phi}. Keeping the derivative squares with $k=i,j$ in each plane gives the first bound.

Set $S=\E[s_i s_j W]$, $X=\E[s_j\partial_iW]$ and $Y=\E[s_i\partial_jW]$. By \eqref{eq:cubictest} each $s_i$ belongs to $L^4(f)$, so $s_i s_j\in L^2(f)$ and $s_i s_j W\in L^1(f)$. Cartesian integration by parts gives $X=Y=-C-S$. Expanding $L_{ij}(s_i s_j)$ and integrating by parts in angle gives $X+Y+\Delta=-S$. Thus $S=\Delta-2C$ and $X=Y=C-\Delta$.

Write $z=\Delta-C$. Projection onto $s_j$ and $s_i$, respectively, gives $\E[(\partial_i W)^2]\ge z^2/I_j$ and $\E[(\partial_j W)^2]\ge z^2/I_i$. The square $A\E[(H_{ij}-wW/A)^2]\ge0$ gives $-A\E[H_{ij}^2]\le-2wC+w^2\tau/A$. Now use $z+C=\Delta$ and $-\alpha z^2+2(w-e)z\le(w-e)^2/\alpha$ to obtain \eqref{eq:pairdual}. This includes $\tau=0$, in which case Cauchy--Schwarz gives $C=0$.
\end{proof}

Write $t_i=T_i=\E v_i^2$ for the coordinate second moments. Then $\sum_i t_i=d$ and $M_i=(d-t_i)/q$. Put
\[
 \omega_{ij}=\frac{M_iM_j}{M_it_i+M_jt_j},\qquad c_i=\frac{7d-(4d+3)M_i}{q}.
\]

\begin{proposition}[Fisher envelope bound]\label{prop:dge4env}
For $d\ge2$ and $f\in\cA_d$,
\begin{equation}\label{eq:dge4env}
\begin{aligned}
 \Phi\le{}&\sum_i(c_iI_i-2M_i^2I_i^2)-\frac2{q^2}\sum_{i<j}\tau_{ij}^2
 -\frac{d-2}{q^2}\sum_{i<j}\tau_{ij}\\
 &-\frac4q\sum_{i<j}\omega_{ij}(1+\tau_{ij})^2.
\end{aligned}
\end{equation}
In particular, $\Phi\le\sum_i(c_iI_i-2M_i^2I_i^2)-d$.
\end{proposition}

\begin{proof}
For $i<j$ put $X=\partial_i w_{ij}-\partial_j u/2$ and $Y=\partial_j w_{ij}+\partial_i u/2$. Since $\E L_{ij}w_{ij}=-\tau_{ij}$ and $\E(v_k\partial_k u)=-1$, we have $\E(Xv_j-Yv_i)=1+\tau_{ij}$. Put $\lambda=\omega_{ij}(1+\tau_{ij})$. Expanding $M_i\E[(X-\lambda v_j/M_i)^2]+M_j\E[(Y+\lambda v_i/M_j)^2]\ge0$ gives
\[
 \omega_{ij}(1+\tau_{ij})^2\le M_i\E X^2+M_j\E Y^2.
\]
By the mixed moment identity in \cref{lem:pairbudget}, the right-hand side equals
\[
 M_i\E[(\partial_iw_{ij})^2]+M_j\E[(\partial_jw_{ij})^2]+
 \delta_{ij}(\Delta_{ij}-C_{ij})+\tfrac14(M_iI_j+M_jI_i).
\]
In \eqref{eq:phi} retain the two indicated derivative squares from each plane and use $\E H_{ii}^2\ge I_i^2$. Summing the linear terms changes $b_i$ to $b_i+(d-M_i)/q=c_i$.

Put $D_{pr}=L_pL_ru$ for coordinate planes $p,r$. Pairing $D_{pr}$ with $D_{rp}$ gives
\[
 \sum_p D_{pp}^2+\tfrac14\sum_{p,r}(D_{pr}-D_{rp})^2\le\sum_{p,r}D_{pr}^2.
\]
The Hessian terms in $D_{pr}-D_{rp}$ cancel by symmetry. Disjoint planes commute. Two planes sharing the index $k$ give $\pm L_{ij}u$, where $i,j$ are their other indices. For fixed $i<j$ there are $d-2$ choices of $k$ and two orders of the planes, so with $W_{ij}=v_i\partial_j u-v_j\partial_i u$,
\[
 \sum_{p,r}(D_{pr}-D_{rp})^2=(d-2)\sum_{i,j}W_{ij}^2=2(d-2)\sum_p(L_pu)^2.
\]
Taking expectations, using $\E D_{pp}^2\ge(\E D_{pp})^2=\tau_p^2$ and multiplying by $2/q^2$ gives the two angular terms in \eqref{eq:dge4env}.

Finally $\sum_{i<j}\omega_{ij}^{-1}=\sum_j(d-t_j)/M_j=dq$. Summing $4\le4\omega_{ij}+\omega_{ij}^{-1}$ over the $dq/2$ pairs gives $\sum_{i<j}\omega_{ij}\ge dq/4$. The last bound follows from $\tau_{ij}\ge0$.
\end{proof}

\section{Reduction to a static inequality}\label{sec:bridge}

We use the identities of \cref{sec:identity} to express the change in Fisher energy along a Gaussian-regularised solution as a time integral of $\Phi$. The sign condition \eqref{eq:static} on $\cA_d$ then implies the monotonicity in \cref{thm:main}. The construction of the weak solution and the integral identity hold independently of this sign condition.

\subsection{Construction and uniqueness of the weak solution}
\begin{proposition}\label{prop:evolution}
Every admissible datum has a weak solution in the sense of \cref{def:weak}. Every such solution conserves mass, momentum and total second moment. Solutions with almost-everywhere equal data agree almost everywhere in velocity at each $t\ge0$. For positive mass, each slice at $t>0$ has a strictly positive representative with a classical score.
\end{proposition}
\begin{proof}
\emph{A discrete approximation.}
Start with the normalised datum and let $f_k$ denote the density before step $k$. Compute $T_k=\int v\otimes v\, f_k$ and $M_k=(d\Id-T_k)/q$. Choose a coordinate plane uniformly and rotate by either sign of $\sqrt{dh}$ with equal probability. Contract by $e^{-h}$ and add an independent centred Gaussian of covariance $(1-e^{-2h})M_k$. The trace stays $d$ and the off-diagonal moments stay zero. In each diagonal entry the recurrence is
\[
 T_{k+1,i}-1=\alpha_h(T_{k,i}-1),\qquad
 \alpha_h=e^{-2h}\bigl(1-2\sin^2(\sqrt{dh})/q\bigr)-(1-e^{-2h})/q.
\]
For small $h>0$, $0<\alpha_h\le1$, whence $\delta\Id\le M_k\le(d/q)\Id$ with $\delta=\min(1,\min_iM_{0,i})>0$.

On tests the step acts as $R_h U_{h,M_k}$, where $R_h$ averages the signed rotations and $U_{h,M}\varphi(v)=\E\varphi(e^{-h}v+Y_h)$ with $Y_h\sim\mathcal N(0,(1-e^{-2h})M)$. For $\varphi\in C_c^\infty$, Taylor expansion along each rotation through the two opposite angles gives
\[
 h^{-1}(R_h-1)\varphi=q^{-1}\sum_{i<j}L_{ij}^2\varphi+O_\varphi(h)
\]
in the supremum norm. The remainder is bounded by the fourth angular derivatives on the ball containing the support. Put $\mathcal U_M\varphi=M:\nabla^2\varphi-v\cdot\nabla\varphi$. Differentiating the Gaussian integral gives $U_{h,M}\varphi-\varphi=\int_0^h U_{r,M}\mathcal U_M\varphi\,\dd r$. The functions $\mathcal U_M\varphi$ have common compact support and uniformly bounded first derivatives for $\delta\Id\le M\le(d/q)\Id$. On a fixed larger ball the contraction tends uniformly to the identity and $\E|Y_r|=O(\sqrt r)$. Outside it, reaching the support requires a Gaussian displacement bounded away from zero. Thus $h^{-1}(U_{h,M}-1)\varphi\to \mathcal U_M\varphi$ uniformly in $M$ and $v$. The rotation average also tends uniformly to the identity on this family. Consequently
\[
 \frac{\langle f_{k+1}-f_k,\varphi\rangle}{h}
 =\biggl\langle f_k,\mathcal U_{M_k}\varphi+q^{-1}\sum_{i<j}L_{ij}^2\varphi\biggr\rangle+o(1)
\]
uniformly in $k$.

\emph{Passage to a weak solution.}
Conditionally on the rotations, the step-$k$ velocity is $e^{-kh}O_kV_0+Y$, with $O_k$ orthogonal. The covariance of $Y$ is a sum of rotated noise covariances with weights $(1-e^{-2h})e^{-2jh}$, $0\le j<k$. It lies between $\delta(1-e^{-2kh})\Id$ and $(d/q)\Id$. In particular $\E|Y|^4\le C_d$. The larger of the two summands in norm is at least half the norm of their sum, so the second-moment tail above $R$ is bounded by four times the initial tail above $R/2$ plus $C_d/R^2$.

Interpolate the laws linearly in time. The one-step estimate gives a common Lipschitz bound for each smooth compact test. Apply Arzel\`a--Ascoli to a countable dense family of such tests and diagonalise over bounded time intervals. Tightness gives a probability measure $\mu_t$ at each time and convergence of the pairings uniformly on compact time intervals. The preceding grid-point law has the same limit, by the Lipschitz bound. Truncating $v_iv_j$ and using the uniform tail bound therefore gives convergence of its second moments to those of $\mu_t$, locally uniformly in $t$. In particular $M_k$ tends to the matrix formed from the limit itself. Summing the one-step formula against a time test and using the normalised generator in \eqref{eq:linear} gives the integral identity of \cref{def:weak} for the measures $\mu_t$.

Fix $t>0$. Both grid points bracketing $t$ have $kh\ge t/2$ for small $h$, so the conditional covariance bound splits off the same Gaussian of covariance $b\Id$, $b=\delta(1-e^{-t/2})/2$, from both laws. Their convex combination is again a convolution with this Gaussian. The residual probabilities have bounded second moment and hence a weakly convergent subsequence. Thus $\mu_t$ has a strictly positive smooth convolution density. To choose the densities jointly in $(t,v)$, fix a probability measure $\nu$ with strictly positive Gaussian density $\gamma$ on $\R^d$. Weak continuity makes $t\mapsto\mu_t$ a Borel probability kernel. Each $\mu_t$ is absolutely continuous with respect to $\nu$, by the convolution representation for $t>0$ and by the prescribed initial density at $t=0$. The Radon--Nikodym theorem for finite kernels on standard Borel spaces gives a jointly measurable $h(t,v)$ such that $\mu_t=h(t,\cdot)\nu$ for each $t$. Choosing $h=0$ where it is infinite changes it only on a $\nu$-null set in each slice. Thus $f_t(v)=h(t,v)\gamma(v)$ is jointly measurable, represents $\mu_t$ for every $t\ge0$, and agrees almost everywhere with $f_{\rm in}$ at $t=0$.

\emph{Conservation and uniqueness.}
For any weak solution choose a smooth $0\le\chi\le1$ of compact support equal to one on $B_1$. Since $1\le\chi(v)+|v|^2$, continuity of $\langle f_t,\chi\rangle$ and the local second-moment bound give a uniform mass bound on $[0,T]$. Thus $B_\varphi(f_t)$ is integrable up to zero for every smooth compact test. By \cref{def:weak} the continuous pairing differs from its time primitive by a constant, fixed by the initial trace. Apply this identity to $1$, $v_i$ and $|v|^2$ times $\chi(v/R)$. The collision kernels are bounded by $C(1+|v|^2+|w|^2)$ uniformly for $R\ge1$, so dominated convergence applies first in $(v,w)$ and then in time. The limiting kernels for $1$ and $v_i$ are zero and $-2(v_i-w_i)$, respectively. For $|v|^2$ the kernel is $2|v-w|^2-4(v-w)\cdot v$ and its symmetrisation in $(v,w)$ vanishes. This proves conservation. At zero initial mass every slice must vanish almost everywhere.

The same cutoffs permit the tests $v_iv_j$. In centred coordinates of unit mass and total second moment $d$, $\sum_{a<b}L_{ab}^2(v_iv_j)=2\delta_{ij}|v|^2-2d v_iv_j$. Thus $\dot T=2M-2T+(2d/q)(\Id-T)$, giving \eqref{eq:Tode}. Two normalised solutions with the same datum now have the same $M(t)$.

For each phase $\theta$, let $z(t,\xi)$ be the difference of their pairings with $\cos(e^t v\cdot\xi+\theta)$. The explicit time derivative of this test cancels $-v\cdot\nabla_v$. Angular derivatives in $v$ transfer to angular derivatives in $\xi$ with opposite sign, hence with the same sign after squaring. These tests follow from the compact-test identity by spatial cutoff and the time product rule: on bounded time and frequency sets the terms are dominated in velocity by $C(1+|v|^2)$. The resulting integral equation makes $z$ absolutely continuous in time, with
\[
 \partial_t z=q^{-1}\sum_{i<j}(L_{ij}^{\xi})^2z-e^{2t}(\xi^{\mathsf T}M(t)\xi)z .
\]
The second moment bounds its first two frequency derivatives locally. For any smooth, nonnegative compact radial weight $\zeta$, angular integration by parts therefore gives
\[
 \tfrac12\frac{\dd}{\dd t}\int\zeta z^2
 =-\int\zeta\biggl(q^{-1}\sum_{i<j}|L_{ij}^{\xi}z|^2+e^{2t}(\xi^{\mathsf T}M\xi)z^2\biggr)\le0
\]
almost everywhere. The initial integral is zero. Exhausting frequency space and using continuity gives $z(t,\xi)=0$ for all $t,\xi$. The phases $0$ and $\pi/2$ determine both parts of the characteristic function, so the laws agree at each time. Inverting the normalisation of \cref{sec:prelim} proves the assertions at positive mass. The zero function supplies existence at zero mass.
\end{proof}
\subsection{Comparison after Gaussian convolution}
Let $f$ be a weak solution in the principal coordinates of \eqref{eq:normal}, $G_\varepsilon$ the centred Gaussian of covariance $\varepsilon\Id$, and $S_r h(v)=r^d h(rv)$. For $\varepsilon>0$ put
\[
 f_t^\varepsilon=f_t*G_\varepsilon,\qquad r=\sqrt{1+\varepsilon},\qquad
 g_t^\varepsilon=S_r f_t^\varepsilon=(S_r f_t)*G_{\varepsilon/(1+\varepsilon)}.
\]

Write $A_M=\sum_iM_i\partial_{ii}+q^{-1}\sum_{i<j}L_{ij}^2$, $\Gamma(h,k)=\nabla h^{\mathsf T}(a+M)\nabla k$ and $\Gamma_2(u)=\tfrac12 A_M\Gamma(u,u)-\Gamma(u,A_Mu)$, holding $M$ fixed in the spatial derivatives. The forward operator in \eqref{eq:linear} is $Q_M h=A_M h+\nabla\cdot(vh)$. For a smooth positive solution with differentiable $M_i(t)$, put $u=\log f$, $\psi=\Gamma(u,u)$ and $e=f\psi$. Then $\partial_t u=A_Mu+\psi+v\cdot\nabla u+d$. Applying the product rule to $e$, the terms $2\Gamma(u,\psi)$ cancel. The drift commutators $[\partial_i,v\cdot\nabla]=\partial_i$ and $[L_{ij},v\cdot\nabla]=0$ give
\begin{equation}\label{eq:local-fisher-balance}
 \partial_t e=Q_M e+\mathcal B,\qquad
 \mathcal B=f\left[\sum_i(\dot M_i+2M_i)(\partial_i u)^2-2\Gamma_2(u)\right].
\end{equation}
The spatial expansion holds pointwise for every smooth $u$. Put $s_i=\partial_i u$, $H_{ij}=\partial_{ij}u$ and $w_{ij}=L_{ij}u$. Give the fields $X=\partial_i,L_p$ weights $\lambda_X=M_i,1/q$, respectively, where $p=(i,j)$ ranges over $i<j$. The product rule gives
\[
 \Gamma_2(u)=\sum_{X,Y}\lambda_X\lambda_Y(XYu)^2+\sum_Y\lambda_Y(Yu)[A_M,Y]u .
\]
Here $[\partial_k,L_{ij}]=\delta_{ki}\partial_j-\delta_{kj}\partial_i$, so $[\sum_k M_k\partial_{kk},L_{ij}]=2(M_i-M_j)\partial_{ij}$. Rotation invariance gives $[\sum_pL_p^2,L_{ij}]=0$. Substitution gives
\[
\begin{aligned}
 \Gamma_2(u)={}&\sum_{i,j}M_iM_jH_{ij}^2+\frac2q\sum_{k,i<j}M_k(\partial_kw_{ij})^2+\frac1{q^2}\sum_{p,r}(L_pL_ru)^2\\
 &+\frac1q\sum_{i<j}(M_i-M_j)\bigl(s_j^2-s_i^2-2s_j\partial_iw_{ij}-2s_i\partial_jw_{ij}+2w_{ij}H_{ij}\bigr).
\end{aligned}
\]
Now suppose the slice belongs to $\cA_d$ and satisfies \eqref{eq:normal}. Put $\Delta_{ij}=I_i-I_j$ and $C_{ij}=\E[w_{ij}H_{ij}]$. The integrations by parts in \cref{lem:pairbudget} give $\E[s_j\partial_iw_{ij}]=\E[s_i\partial_jw_{ij}]=C_{ij}-\Delta_{ij}$, so the last parenthesis has mean $3\Delta_{ij}-2C_{ij}$. In \eqref{eq:local-fisher-balance}, keep $-4q^{-1}\sum_{i<j}(M_i-M_j)(\Delta_{ij}-C_{ij})$. The remaining pair sum is $-2q^{-1}\sum_{i<j}(M_i-M_j)\Delta_{ij}=-2dq^{-1}\sum_i(M_i-1)I_i$. Hence the coefficient of $I_i$ is $\dot M_i+2M_i-2d(M_i-1)/q=b_i$ by \eqref{eq:Tode}, giving $\int\mathcal B=\Phi(f_t)$ as in \eqref{eq:phi}.

For the derivative of the full integral at $t>0$, assume additionally that $e_t\in L^1$. For $X=\partial_i$ and $b=\partial_i\psi,\psi\partial_i u,v_i\psi$, require $b,Xb,bXu\in L^1(f\,\dd v)$. For $X=L_{ij}$ and $b=X\psi,\psi Xu$, require $Xb,bXu\in L^1(f\,\dd v)$. Each $X(fb)=f(Xb+bXu)$ has integral zero by testing against the radial cutoffs $\chi_R$ of \cref{lem:closedD}. Indeed $|\partial_i\chi_R|\le C/R$ and $L_{ij}\chi_R=0$. Expanding $Q_M e$ into these derivatives gives $\int Q_M e=0$. If $|\partial_\tau e_\tau(v)|\le B(v)$ for almost every $v$ and all $\tau$ near $t$, with $B\in L^1(\dd v)$, differentiation under the integral gives $(\dd/\dd t)\int e=\Phi(f_t)$.

\begin{proposition}\label{prop:gaussian-comparison}
Under \eqref{eq:static}, $\cF(f_t^\varepsilon)\le\cF(f_s^\varepsilon)$ for every $\varepsilon>0$ and $0\le s\le t$.
\end{proposition}
\begin{proof}
Radial convolution commutes with $L_{ij}$, and $(v_i f)*G_\varepsilon=v_i f^\varepsilon+\varepsilon\partial_i f^\varepsilon$. Hence $f^\varepsilon$ solves $\partial_t f^\varepsilon=A_{M+\varepsilon\Id}f^\varepsilon+\nabla\cdot(vf^\varepsilon)$ in distributions. Dilation commutes with rotations and the drift and divides the translation coefficients by $r^2$. Thus $g=g^\varepsilon$ solves the normalised equation at the same times, with second moments $(T_i+\varepsilon)/(1+\varepsilon)$ and its own coefficients $\widetilde M_i=(M_i+\varepsilon)/(1+\varepsilon)$. In particular $\dot{\widetilde M}_i=4d(1-\widetilde M_i)/q$.

For the time regularity, pair $f_t(y)$ with translates of a polynomial times a Gaussian. On compact sets of translation parameters these kernels and their derivatives decay rapidly in $y$, uniformly. Cutoffs, unit mass and compact-test continuity give joint continuity of the pairings up to $t=0$. The weak equation differentiates each pairing at $t>0$ by applying $A_M-y\cdot\nabla_y$ to the kernel. This operator and spatial derivatives give finite sums of the same type of kernel, with smooth coefficients by \eqref{eq:Tode}. Iteration gives joint smoothness of $g$ for positive times and the forward equation pointwise.

At any $t\ge0$, $g$ is the density of $V=Y+Z$, where $Y$ has density $S_rf_t$ and $Z$ is an independent Gaussian of covariance $\eta\Id$, $\eta=\varepsilon/(1+\varepsilon)$. Set $R_\alpha=\partial^\alpha g/g$ and $P_\alpha=\partial^\alpha G_\eta/G_\eta$. Since $R_\alpha(V)=\E[P_\alpha(Z)\mid V]$, conditional Jensen gives, for $m=2,4$,
\[
 \E[(1+|V|^2)|R_\alpha(V)|^m]\le\E[(1+|V|^2)|P_\alpha(Z)|^m]\le C_{d,\eta,\alpha,m}(1+\E|Y|^2).
\]
The last bound uses independence, $|V|^2\le2|Y|^2+2|Z|^2$ and the Gaussian moments of the polynomial $P_\alpha$. For $u=\log g$, $s_i=R_i$ and $H_{ij}=R_{ij}-R_iR_j$, so
\[
 \int(1+|v|^2)g\bigl(|\nabla u|^2+|\nabla u|^4+|\nabla^2u|^2\bigr)\,\dd v\le C_{d,\eta}(1+\E|Y|^2).
\]
In particular $w_{ij}=v_i s_j-v_j s_i$ and $\partial_kw_{ij}=\delta_{ki}s_j-\delta_{kj}s_i+v_iH_{kj}-v_jH_{ki}$ belong to $L^2(g\,\dd v)$ uniformly in time at fixed variance. Since $X\sqrt g=\tfrac12\sqrt g\,Xu$ for each field $X$, the graph formula gives $\int e_t=\cF(g_t)<\infty$ at every $t\ge0$, with $e=g\Gamma(u,u)$ and $M=\widetilde M$.

Set $U=2q^{-2}g\sum_{p,r}|L_pL_ru|^2$, with $p,r$ independently ranging over the coordinate planes, and $\widetilde{\mathcal B}=\mathcal B+U$. The pointwise expansion above, the weighted bounds and Cauchy--Schwarz give $e,\widetilde{\mathcal B}\in L^1([0,t]\times\R^d)$ for each $t<\infty$. Testing \eqref{eq:local-fisher-balance} against $\chi_R$ for $0<s<t$ gives
\[
 \int\chi_R e_t+\int_s^t\!\int\chi_R U
 =\int\chi_R e_s+\int_s^t\!\int(\chi_R \widetilde{\mathcal B}+eD_R),\qquad
 D_R=\sum_i\widetilde M_i\partial_{ii}\chi_R-v\cdot\nabla\chi_R .
\]
Here $L_{ij}\chi_R=0$ and $D_R$ is uniformly bounded for $R\ge1$ and tends pointwise to zero. At fixed $R$, positivity of $g$ and continuity of its spatial derivatives up to time zero give $\int\chi_R e_s\to\int\chi_R e_0$. Monotone convergence of the time integral of $\chi_R U$ and dominated convergence on the right extend the identity to $s=0$.

Dropping $\int\chi_R e_t\ge0$ gives a bound on $\int_0^t\!\int\chi_R U$ independent of $R$. Fatou's lemma and Fubini give $U\in L^1([0,t]\times\R^d)$ and $L_pL_ru\in L^2(g_\tau\,\dd v)$ for almost every interior time $\tau$. Thus $g_\tau\in\cA_d$ at these times and $\int\mathcal B_\tau=\Phi(g_\tau)$. Dominated convergence in the cutoff identity and the graph formula at both endpoints now give, for every $0\le s\le t$,
\[
 \cF(g_t)-\cF(g_s)=\int_s^t\Phi(g_\tau)\,\dd\tau .
\]
Assuming \eqref{eq:static} on $\cA_d$, the integrand is nonpositive for almost every $\tau$, so the integral is nonpositive. Velocity dilation preserves $\cF$, proving the comparison.
\end{proof}

\subsection{Removing the convolution}
\begin{proposition}\label{prop:graph-comparison}
Under \eqref{eq:static}, every weak solution satisfies $\cF(f_t)\le\cF(f_s)$ for $0\le s\le t$.
\end{proposition}
\begin{proof}
Assume \eqref{eq:normal} and take $\varepsilon_n=1/(n+1)$. At either endpoint, convolution converges in $L^1$, preserves mass and momentum, and adds $d\varepsilon_n$ to the total second moment. Thus \cref{lem:closedD} gives $\cF(f_t)\le\liminf_n\cF(f_t^{\varepsilon_n})$ and likewise at $s$. If $\cF(f_s)<\infty$, put $\rho=\sqrt{f_s}$. For $X=\partial_i$ or $L_{ij}$, the weak product rule gives $Xf_s=2\rho X\rho\in L^1$. Commuting $X$ with radial convolution and applying Cauchy--Schwarz gives pointwise
\[
 |X\sqrt{f_s^\varepsilon}|^2=\frac{|(\rho X\rho)*G_\varepsilon|^2}{f_s^\varepsilon}\le (|X\rho|^2)*G_\varepsilon .
\]
Integration contracts each derivative norm. The covariance changes from $T$ to $T+\varepsilon\Id$, hence $M_{f_s^\varepsilon}=M_{f_s}+\varepsilon\Id$. The graph formula therefore gives
\[
 \cF(f_s^\varepsilon)\le\cF(f_s)+4\varepsilon\|\nabla\rho\|_2^2.
\]
Thus the energy is recovered at $s$, and \cref{prop:gaussian-comparison} gives
\[
 \cF(f_t)\le\liminf_n\cF(f_t^{\varepsilon_n})\le\limsup_n\cF(f_s^{\varepsilon_n})=\cF(f_s).
\]
An infinite starting value gives the comparison immediately. For mass $m>0$ the original energy is $m^2$ times the normalised energy at time $mt$. At zero mass it vanishes.
\end{proof}

\subsection{The initial time}
\begin{proposition}\label{prop:initial}
Every weak solution satisfies $\cF(f_{\rm in})\le\liminf_{t\downarrow0}\cF(f_t)$. Under \eqref{eq:static}, $D_{\rm cl}(f_t)\to D_{\rm cl}(f_{\rm in})$ in $[0,\infty]$ as $t\downarrow0$.
\end{proposition}
\begin{proof}
At positive mass use the normalisation \eqref{eq:normal}. Given a sequence $t_n>0$ tending to zero with finite energy liminf, pass to a bounded-energy subsequence realising it. By \eqref{eq:Tode}, $M(t_n)\ge\delta\Id$ with $\delta=\min\{1,\min_i M_i(0)\}>0$, so $4\delta\|\nabla\sqrt{f_{t_n}}\|_2^2\le\cF(f_{t_n})$. For $\rho_n=\sqrt{f_{t_n}}$, the translation bound $\|\rho_n(\cdot-z)-\rho_n\|_2\le |z|\|\nabla\rho_n\|_2$ and $\|\rho_n\|_2=1$ give $\|f_{t_n}(\cdot-z)-f_{t_n}\|_1\le2|z|\|\nabla\rho_n\|_2$. Averaging over $G_\varepsilon$ gives
\[
 \|f_{t_n}^\varepsilon-f_{t_n}\|_1\le2\sqrt{d\varepsilon}\,\|\nabla\sqrt{f_{t_n}}\|_2.
\]
For fixed $v,\varepsilon$, the test $G_\varepsilon(v-\cdot)$ is a uniform limit of compact tests. The common unit mass and compact-test continuity therefore give $f_{t_n}^\varepsilon(v)\to(f_{\rm in}*G_\varepsilon)(v)$. These nonnegative densities all have integral one, so the convergence is in $L^1$ by Scheff\'e's lemma. In the triangle inequality, first let $n\to\infty$ for fixed $\varepsilon$, then use the displayed uniform estimate and the $L^1$ approximate-identity property of $G_\varepsilon$ on $f_{\rm in}$. This gives $f_{t_n}\to f_{\rm in}$ in $L^1$. Conservation and \cref{lem:closedD} give the lower bound. Assuming \eqref{eq:static} on $\cA_d$, apply \cref{prop:graph-comparison} for the upper bound at $s=0$ and use $f_0=f_{\rm in}$ almost everywhere. Subtraction of $dm^2$ gives the limit, including $+\infty$. At zero mass both quantities vanish.
\end{proof}

\subsection{Conclusion of the reduction}
\begin{theorem}\label{thm:reduction}
Fix $d\ge2$. If \eqref{eq:static} holds on $\cA_d$, the conclusions of \cref{thm:main} hold in dimension $d$.
\end{theorem}
\begin{proof}
Conservation of mass and \cref{prop:graph-comparison,prop:initial} give \eqref{eq:main}. Fix $t>0$ in principal coordinates and choose $0<\delta\le\min(1,\min_i M_i(0))$. By \eqref{eq:Tode}, $(M_i+\varepsilon)/(1+\varepsilon)\ge\delta$ for every time and variance.

For $g^\varepsilon$ form the vector $z$ with entries $(\widetilde M_i I_i)_{1\le i\le d}$ and $(\tau_{ij}/q)_{i<j}$. Their sum is $Y=\cF(g^\varepsilon)$. Put $t_0=t/2$. At almost every $\tau\in[t_0,t]$, dropping the last two nonpositive terms of \eqref{eq:dge4env} gives
\[
 \Phi(g_\tau^\varepsilon)\le\sum_i c_i I_i-2\sum_\ell z_\ell^2,\qquad
 c_i=\frac{7d-(4d+3)\widetilde M_i}{q}\le\frac{7d}{q}\le A\widetilde M_i,\qquad A=\frac{7d}{q\delta}.
\]
All $N=d+\binom d2=d(d+1)/2$ entries are nonnegative. Thus $\sum_i c_i I_i\le AY$ and $\sum_\ell z_\ell^2\ge Y^2/N$, which give
\[
 \Phi(g_\tau^\varepsilon)\le A Y(\tau)-c_d Y(\tau)^2,\qquad c_d=2/N=\frac4{d(d+1)}.
\]
The balance on every subinterval of $[t_0,t]$ makes $Y$ absolutely continuous with $Y'=\Phi(g^\varepsilon)$ almost everywhere. For $\eta,c>0$ set $B(x)=(A+\eta)/c_d+1/[c_d(x-t_0+c)]$, choosing $c$ small enough that $B(t_0)>Y(t_0)$. Then $B'>AB-c_d B^2$, so $(Y-B)'\le[A-c_d(Y+B)](Y-B)$ almost everywhere. An integrating factor gives $Y(t)\le B(t)\le(A+\eta)/c_d+2/(c_d t)$. Let $\eta\downarrow0$. Dilation invariance and lower semicontinuity pass the bound $A/c_d+2/(c_d t)$ to $\cF(f_t)$ as $\varepsilon\downarrow0$.

Returning to the original variables preserves finiteness. For positive mass, \cref{prop:evolution,lem:closedD} give a positive representative with a classical score at every $t>0$, and $D(h,p)=D_{\rm cl}(f_t)<\infty$ for every such pair. Thus \eqref{eq:main} compares the classical integrals for every pair of choices at $0<s\le t$.
\end{proof}

\section{The static inequality in dimension three}\label{sec:d3}

In dimension three we optimise the bound in each coordinate plane while retaining its dependence on the angular energy $\tau_{ij}$. The spherical estimate then controls the sum of the angular contributions.

Here $q=2$, $M_i=M(T_i)$ with $M(t)=(3-t)/2$, $b_i=9-7M_i$ and $B=B_3$. By \cref{lem:sphere}, $\cK\ge J+J^2$. In each plane write $\delta=\delta_{ij}=e$, with $A,\alpha,\Delta$ as in \cref{lem:pairbudget}. The quadratic in \eqref{eq:pairdual} attains its minimum at $w=A(\alpha\Delta+\delta)/(A+\alpha\tau)$, giving
\[
 R_{ij}(\tau)=\frac{\delta^2\tau-2A\delta\Delta-\alpha A\Delta^2}{A+\alpha\tau}\qquad(\tau\ge0).
\]
Summing the pair bounds yields
\begin{equation}\label{eq:d3U}
 \Phi\le B+\sum_{i<j}R_{ij}(\tau_{ij})-J-J^2 .
\end{equation}

\subsection{The price envelope}
For each plane put $\eta_{ij}=\delta^2/A$. With $y_+=\max(y,0)$ and $\psi_c(y)=y^2-(|y|-c)_+^2$ for $c\ge0$, completion of the square gives
\begin{equation}\label{eq:Fexact}
 F_{ij}=\min_{|w|\le|\delta|}\left\{\frac{(w-\delta)^2}{\alpha}-2w\Delta\right\}
       =\frac{\delta^2-\psi_{|\delta|}(\alpha\Delta+\delta)}{\alpha} .
\end{equation}
The minimiser $w_*$ is the projection of $\alpha\Delta+\delta$ onto $[-|\delta|,|\delta|]$. For equal weights, $w_*=F_{ij}=0$.

Testing the preceding unrestricted minimum at $w_*$ gives $R_{ij}(\tau)\le F_{ij}+\eta_{ij}\tau$. Set $\eta=\max_{i<j}\eta_{ij}$ and $P(T,I)=B+\sum_{i<j}F_{ij}+\eta^2$. Since $\sum_{i<j}\tau_{ij}=2J$, \eqref{eq:d3U} becomes
\begin{equation}\label{eq:d3P}
 \Phi\le P(T,I)-J-(J-\eta)^2 .
\end{equation}
We shall bound $P$ for all $T_i>0$, $\sum_iT_i=3$ and $I_i\ge1/T_i$.

\subsection{Moderate temperatures}
\begin{lemma}\label{lem:d3moderate}
If $\max_i T_i\le3/2$, then $P(T,I)\le-D_T/6$, where $D_T=\sum_i(T_i-1)^2$.
\end{lemma}
\begin{proof}
Taking $w=\delta$ in \eqref{eq:Fexact} gives $F_{ij}\le-2(M_i-M_j)(I_i-I_j)$. Since $\sum_i M_i=3$, summing yields $P\le\sum_i[(15-13M_i)I_i-2M_i^2I_i^2]+\eta^2$. Write $r_i=1/T_i$, $S=\sum_i r_i$ and $E_T=\sum_i(r_i-1)^2$. Since $(3-T_i)(3-T_j)\ge T_iT_j$ and $|T_i-T_j|\le3$,
\[
 \eta_{ij}=\frac{(T_i-T_j)^2}{4(3-T_i)(3-T_j)}\le\tfrac34|r_i-r_j|,
 \qquad (r_i-r_j)^2\le2E_T .
\]
Thus $\eta^2\le9 E_T/8$. At $I_i=r_i$ the $i$th summand is $-\tfrac92(r_i-1)^2-\tfrac{21}2(r_i-1)$.

Increasing $I_i$ by $v\ge0$ from $r_i$ adds $g(T_i)v-2M_i^2v^2$, where $g(t)=(11t+3)/2-9/t$. The maximal gain is $g(T_i)_+^2/(8M_i^2)$. It is zero for $T_i\le1$, since $g$ is increasing and $g(1)=-2$. For $1\le t\le3/2$, the identity $\tfrac{15}2(t-1)-g(t)=(3-2t)(3-t)/t\ge0$ and $M(t)\ge3/4$ bound the gain by $\tfrac{25}2(t-1)^2$.

Let $P_T=\sum_{T_i>1}(T_i-1)^2$. The positive deviations sum to the absolute sum of the negative ones, of which there are at most two. Cauchy--Schwarz gives $P_T\le2(D_T-P_T)$, or $P_T\le2D_T/3$. Also $E_T\ge4D_T/9$ and $S-3=\sum_i(T_i-1)^2/T_i\ge2D_T/3$. Combining the bounds gives
\[
 P\le-\tfrac{27}8 E_T-\tfrac{21}2(S-3)+\tfrac{25}2P_T\le-D_T/6 .\qedhere
\]
\end{proof}

\subsection{Two saturated cold directions}
If a temperature exceeds $3/2$, permute the coordinates to write $T=(t,a,b)$ with $t>3/2$ and $0<a\le b$. Let $s=(a+b)/2=(3-t)/2$, $h=3/s-3$, and abbreviate $P(T,(x,y,z))$ by $P(x,y,z)$.

\begin{lemma}\label{lem:d3twoCR}
For $x\ge1/t$, $P(x,1/a,1/b)\le-15h^2/64-(a-b)^2/(ab)$.
\end{lemma}
\begin{proof}
Write $d_\theta(y)=(9-7M(\theta))y-2M(\theta)^2y^2$. The choice $w=\delta$ bounds the cold--cold pair by $-(a-b)^2/(ab)$. For $\theta=a,b$, the choices $w=\delta,-\delta$ for the pair $(t,\theta)$ bound its contribution plus $d_\theta(1/\theta)$ by $C_+(\theta),C_-(\theta)$ respectively, where $k=(3-t)x-1$, $D_\theta=3+k\theta>0$ and
\[
\begin{aligned}
 C_+(\theta)&=4+(t-\theta)x+(3/2-t)/\theta-9/(2\theta^2),\\
 C_-(\theta)&=2-(t-\theta)x+(t+3/2)/\theta-9/(2\theta^2)+x(t-\theta)^2/D_\theta .
\end{aligned}
\]
Let $\mathcal D g=g(a)+g(b)-2g(s)$. The two reciprocal terms in $C_+$ have nonpositive midpoint differences by convexity of $1/\theta$ and $1/\theta^2$ and $t>3/2$. For $C_-$, direct subtraction gives
\[
 \mathcal D C_-=(a-b)^2\left[\frac{t+3/2}{2abs}
 -\frac{9(2s^2+ab)}{4a^2b^2s^2}+\frac{x(kt+3)^2}{2D_aD_bD_s}\right].
\]
Here $kt+3=(3-t)(tx+1)>0$. Since $D_\theta-\theta(kt+3)/t=3(1-\theta/t)\ge0$ for $\theta=a,b$ and $D_s\ge(3-t)sx$, the last fraction is at most $t^2/[2ab(3-t)s]$. With $t=3-2s$ the resulting bracket equals $-(3ab+18s)/(4a^2b^2s)$. Thus $\mathcal D C_\pm\le0$, including at $a=b$.

Put $X=sx$, so $(h+1)X=tx\ge1$, and $B_0=-2X^2+(3h+2)X-h^2-5h$. At the midpoint, $d_t(x)+2C_\pm(s)$ equals, respectively,
\[
 B_0+2h(X-1),\qquad B_0+\frac{2h(h+2-2X^2)}{2X+h+2} .
\]
For $X\le1$, write $u=1-X\ge0$. The first expression becomes $-h^2-2h+(2-5h)u-2u^2\le-h^2/4$ since $h>1$. For $X\ge1$, add $h^2/4$ to the second and multiply by $2X+h+2$. With $v=X-1$ this gives
\[
 -\bigl(4v^3-\tfrac32h^2v+\tfrac34h^3\bigr)-12v^2-(2h+8)v-3h^2-8h\le0 .
\]
For $0\le v\le h/2$ the term $3h^3/4$ dominates $3h^2v/2$. For $v=h/2+y$, $y\ge0$, the parenthesis is $4y^3+6hy^2+3h^2y/2+h^3/2$.

All three weights lie in $[s,3/2]$. For a pair with larger-to-smaller ratio $R\ge1$, its price is $(R+1/R-2)/4$, increasing in $R$. Hence $\eta\le(3/2-s)^2/(6s)\le h/8$, the last inequality reducing to $s(3-4s)\ge0$. Adding the cold--cold term and $\eta^2$ to $-h^2/4$ proves the lemma.
\end{proof}

\subsection{High temperature with general Fisher informations}
\begin{lemma}\label{lem:d3hot}
If $\max_i T_i=t>3/2$, then for all $I_i\ge1/T_i$,
\[
 P(T,I)\le-\frac{15h^2}{64}=-\frac{135}{64}\left(\frac{t-1}{3-t}\right)^2<0 .
\]
\end{lemma}
\begin{proof}
The derivative of $\psi_c$ is twice the projection onto $[-c,c]$, including at the endpoints. Thus \eqref{eq:Fexact} gives $\partial_{I_i}F_{ij}=2M_j(w_*-\delta)^2/(\alpha I_i)^2-2w_*$. Since $\alpha I_i\ge2M_j$, this is bounded above by $(w_*-\delta)^2/(2M_j)-2w_*$. The maximum of this quadratic on the interval occurs at an endpoint, giving
\begin{equation}\label{eq:Fderiv}
 \partial_{I_i}F_{ij}\le\begin{cases}2M_i(M_i-M_j)/M_j,&M_i\ge M_j,\\ 2(M_j-M_i),&M_i\le M_j.\end{cases}
\end{equation}
Write $m_\theta=M(\theta)$. By the mean value theorem, raising the first cold Fisher coordinate from $1/a$ to $1/a+v$, $v\ge0$, changes each of its pair terms by at most the bound in \eqref{eq:Fderiv} times $v$. Adding the change in $d_a$ bounds the change in $P$ by
\[
 \left[9-11m_a+2m_a^2\left(\frac1{m_t}+\frac1{m_b}-\frac2a\right)\right]v-2m_a^2v^2\le0 ,
\]
since $a\le m_t,m_b$ and $a<3/4$. Hence $P(x,y,z)\le P(x,1/a,z)$ for $y\ge1/a$.

For $z=1/b+v$, $v\ge0$, the same estimate bounds the change from $P(x,1/a,1/b)$ by $Lv-2m_b^2v^2$, where
\[
 L=\frac{11b-2a-9}{2}-\frac{(3-b)^2a}{(a+b)b}\le3(b-a).
\]
Indeed $5b+4a\le5(a+b)<15/2<9$, so the first term of $L$ is at most $3(b-a)$. Completing the square bounds the gain by $9(b-a)^2/(8m_b^2)$. Since $t\ge a+b$, we have $3-b=t+a\ge2a+b\ge2\sqrt{2ab}$ and hence $8m_b^2=2(3-b)^2\ge16ab\ge9ab$. The gain is therefore at most $(a-b)^2/(ab)$ and cancels the last term of \cref{lem:d3twoCR}, including when $a=b$.
\end{proof}

\begin{theorem}\label{thm:d3}
For every $f\in\cA_3$, with $T_i=\E[v_i^2]$ in the given coordinates,
\[
 \Phi(f)\le-J-(J-\eta)^2\le0 .
\]
The first inequality is strict if any $T_i>3/2$.
\end{theorem}
\begin{proof}
The actual moments satisfy $T_i>0$, $\sum_iT_i=3$ and $I_i\ge1/T_i$. Apply \cref{lem:d3moderate,lem:d3hot} in \eqref{eq:d3P} and use $J\ge0$.
\end{proof}

\section{The static inequality in dimensions \texorpdfstring{$d\ge4$}{d >= 4}}\label{sec:d4}

Let $d\ge4$ and write $t_i=T_i$. When every $t_i\le d-1$, the Fisher envelope in \cref{prop:dge4env} lets us maximise separately over the Fisher informations $I_i$. When one temperature exceeds $d-1$, we instead use \eqref{eq:pairdual} to retain the dependence on the other temperatures.

For $0<t<d$ write $M(t)=(d-t)/q$ and $c(t)=[7d-(4d+3)M(t)]/q$. The slope of $c(t)I-2M(t)^2I^2$ at $I=1/t$ is $P_d(t)/(q^2t)$, where
\[
 P_d(t)=(4d-1)t^2+d(3d-2)t-4d^2.
\]
This polynomial increases strictly on $[0,\infty)$, with $P_d(1)=-q^2$ and $P_d(2)=2d^2+12d-4>0$. Let $\kappa\in(1,2)$ be its positive root. Write $g_d(t)$ for the maximum over $I\ge1/t$. For $t\le\kappa$ it is the endpoint value $F_d(t)$, and for $t>\kappa$ it is the value $G_d(t)$ at the vertex $c(t)/(4M(t)^2)$.
The two values are
\[
 F_d(t)=\frac{4d+1+3d(d-2)/t-2d^2/t^2}{q^2},\qquad
 G_d(t)=\frac18\Bigl[\frac{7d}{d-t}-\frac{4d+3}{q}\Bigr]^2.
\]
These values agree at $\kappa$ since $8q^2t^2(d-t)^2(G_d-F_d)=P_d(t)^2$. On $(0,\kappa]$ we have $F_d''(t)=6d[(d-2)t-2d]/(q^2t^4)<0$. If $B(t)$ denotes the bracket in $G_d$, then $B(t)\ge3$ for $1\le t<d$ and $B''(t)=14d/(d-t)^3>0$. Hence $4G_d''=(B')^2+BB''>0$ on $[\kappa,d)$.

\begin{proposition}[Unpriced pair bound]\label{prop:dge4pair}
For $f\in\cA_d$,
\[
 \Phi\le\mathcal E(t,I):=\sum_i(b_iI_i-2M_i^2I_i^2)
 +\sum_{i<j}\frac{(M_i-M_j)^2I_iI_j}{q(M_iI_i+M_jI_j)},
\]
where $b_i=[6d-(4d+2)M_i]/q$ as in \eqref{eq:phi}.
\end{proposition}
\begin{proof}
Set $w=0$ in \eqref{eq:pairdual}. The bound $e^2/\alpha$ is the displayed pair term. Sum over the pairs in \cref{lem:pairbudget} and use $\cK\ge0$.
\end{proof}

\begin{theorem}\label{thm:dge4}
For $d\ge4$ and $f\in\cA_d$, $\Phi(f)\le0$. If any $T_i>d-1$, then $\Phi(f)<0$.
\end{theorem}

\begin{proof}
\emph{Ordinary region.} Suppose $\max_i t_i\le q$. Since $I_i\ge1/t_i$, \cref{prop:dge4env} reduces the claim to $\sum_i g_d(t_i)\le d$. If every $t_i\le\kappa$, Jensen's inequality gives $\sum_i F_d(t_i)\le dF_d(1)=d$.

For two hot coordinates $u,v>\kappa$, put $y=u+v-\kappa\le d-\kappa<q$. Both $u,v$ lie on $[\kappa,y]$ and $u+v=\kappa+y$. The sum of the chord values of $G_d$ at $u,v$ is therefore $G_d(\kappa)+G_d(y)$, which bounds $G_d(u)+G_d(v)$ by convexity. Replacing $(u,v)$ by $(\kappa,y)$ preserves the total temperature and removes one hot coordinate. Since $F_d(\kappa)=G_d(\kappa)$, the envelope sum can only increase. Repeat until one hot value $x$ remains.

The other $q$ temperatures have mean $(d-x)/q$. Put $z=q/(d-x)$ and $z_\kappa=q/(d-\kappa)$, so $1<z_\kappa\le z\le q$. Jensen's inequality on these $q$ coordinates gives $\sum_i g_d(t_i)-d\le G_d(x)+qF_d(1/z)-d=p_d(z)$, where
\[
 8q^2p_d(z)=d^2(65-16d)z^2+2d(12d^2-64d+3)z-8d^3+64d^2-8d+1.
\]
To evaluate the left endpoint, substitute $z=q/(d-x)$ in this quadratic. The difference between $8q^2 x(d-x)^2p_d(q/(d-x))$ and $-4d^2(x-1)^2[(6d-13)x+d(5d+2)]$ factors as $P_d(x)[(4d^2+4d-1)x-d(5d+2)]$. Setting $x=\kappa$ and writing $r=d-4\ge0$ gives
\[
\begin{aligned}
 8q^2\kappa(d-\kappa)^2p_d(z_\kappa)&=-4d^2(\kappa-1)^2[(6d-13)\kappa+d(5d+2)]<0,\\
 8q^2p_d(q)&=-16r^5-199r^4-930r^3-2033r^2-2102r-839<0.
\end{aligned}
\]
For integer $d\ge5$ the leading coefficient has the sign of $65-16d<0$, so $p_d'(z)\le p_d'(1)=-d(4d+3)/(4q)<0$ for $z\ge1$. For $d=4$ the coefficient is positive and $p_4(3)=-839/72$ by the second identity. Convexity bounds $p_4$ by the chord between its two negative endpoint values. Thus $p_d<0$ on $[z_\kappa,q]$ in both cases.

\emph{Extreme region.} Label the extreme coordinate $t_0=d-h>q$. The other $q$ temperatures $\theta_j$ are positive and sum to $h$, so $0<h<1$. In \cref{prop:dge4pair} the derivative of a pair term with respect to $I_i$ is bounded by $(M_i-M_j)^2/(qM_j)$ since its remaining factor $[M_jI_j/(M_iI_i+M_jI_j)]^2$ is at most one. Summing and using $I_i\ge1/t_i$ gives, for $a=M_i$,
\[
 \partial_{I_i}\mathcal E\le\frac{7d-(6d+2)a}{q}
 +a^2\left(\sum_j\frac1{d-t_j}-\frac4{t_i}\right).
\]
The chord bound $1/(d-z)\le1/d+z/(dt_i)$ holds on $0\le z\le d-t_i$. Since $\sum_{j\ne i}t_j=d-t_i$, summation gives $\sum_{j\ne i}1/(d-t_j)\le1/t_i+(d-2)/d$. For $t_i\le1$, the bracket is at most $1/(d-t_i)-3/t_i+(d-2)/d\le1/q-2-2/d<0$. With $a\ge1$ this yields $\partial_{I_i}\mathcal E\le(d-2)/q+1/q-2-2/d=-(d+2)/d$. Thus lowering each cold Fisher coordinate to $\widehat I_j=1/\theta_j$ gives $\mathcal E(t,I)\le\mathcal E(t,\widehat I)$, with $\widehat I_0=I_0$.

Put $\sigma=\sum_j1/\theta_j$ and $\sigma_2=\sum_j1/\theta_j^2$. Complete the square in the hot diagonal term. For a hot--cold pair, the denominator is at least $qM_0 I_0$, giving an upper bound $(M_0-M_j)^2/(qM_0\theta_j)$. The contributions to $\mathcal E(t,\widehat I)$ satisfy
\[
\begin{aligned}
 \text{hot diagonal term}&\le\tfrac18[6d/h-(4d+2)/q]^2,\\
 \text{cold diagonal sum}&=[4dq+2d(d-2)\sigma-2d^2\sigma_2]/q^2,\\
 \text{hot--cold pair sum}&\le[(d-h)^2\sigma/h-2q(d-h)/h+1]/q^2.
\end{aligned}
\]
For cold temperatures $a,b\le1$, the pair term is
\[
 \frac{(a-b)^2}{q^2[d(a+b)-2ab]}\le\frac{a+b}{q^3},
\]
since $d(a+b)-2ab\ge q(a+b)$. Half the ordered sum over all cold indices is therefore at most $h/q^2\le1/q^2$.

Cauchy--Schwarz gives $\sigma\ge s_0=q^2/h$ and $\sigma_2\ge\sigma^2/q$. In the numerator over $q^2$, the terms depending on these moments are bounded above by $R(\sigma)$, where
\[
 R(s)=\bigl[2d(d-2)+(d-h)^2/h\bigr]s-2d^2s^2/q.
\]
Here $hR'(s_0)=d^2(1-4q)+2d(d-3)h+h^2\le d^2(7-4d)-6d+1<0$. Thus $R(\sigma)-R(s_0)=(\sigma-s_0)R'(s_0)-2d^2(\sigma-s_0)^2/q\le0$.
The hot diagonal term, the constant parts of the two sums above, and $R(s_0)$ give
\[
 \mathcal E(t,I)\le\frac18\left[\frac{6d}{h}-\frac{4d+2}{q}\right]^2
 +\frac{4dq-2q(d-h)/h+2+R(q^2/h)}{q^2}
 =\frac{\mathsf A}{h^2}+\frac{\mathsf B}{h}+\mathsf C,
\]
where
\[
 \mathsf A=\frac{d^2(15-4d)}2,\qquad \mathsf B=\frac{d(2d^2-14d+1)}q,\qquad
 \mathsf C=\frac{14d^2-4d+3}{2q^2}.
\]
Here $\mathsf A<0$. With $r=d-4$,
\[
\begin{aligned}
 2\mathsf A+\mathsf B&=-\tfrac dq(4r^3+27r^2+53r+35)<0,\\
 2q^2(\mathsf A+\mathsf B+\mathsf C)&=-4r^5-53r^4-274r^3-701r^2-910r-485<0.
\end{aligned}
\]
Thus $\mathsf A z^2+\mathsf B z+\mathsf C$ is negative at $z=1$ and decreases for $z\ge1$. Taking $z=1/h>1$ gives $\Phi\le\mathcal E(t,I)<0$.
\end{proof}

\section{The static inequality in dimension two}\label{sec:d2}

\subsection{Notation and the exact formula}
Assume the static regularity, unit mass and total second moment of \cref{def:staticclass} in dimension two. For the exact formula below the mean may be arbitrary. Put $a=\E v_2^2$, $b=\E v_1^2$, $\Theta=\diag(b,a)$ and $M=\diag(a,b)$. Then $a,b>0$, $a+b=2$, and $\E(v_1v_2)$ may be nonzero. Write $p=\nabla u$, $h_i=H_{ii}$, $L=L_{12}$ and $F=Lu$. In these coordinates $J=\E F^2$. Set
\[
 \nu=\E(FH_{12}),\qquad \Lambda_i=\E[(\partial_i F)^2],\qquad K_F=\E[(LF)^2].
\]
With $|z|_M^2=z^{\mathsf T}Mz$, write $TT=\E\tr(MHMH)$, $S_M=\E(|\nabla F|_M^2+|Lp|_M^2)$ and $\Eord=TT+S_M+K_F$. The planar form is $Q_L=\Eord+6(aI_1+bI_2)-8(I_1+I_2)$.

\begin{proposition}\label{prop:d2formula}
Let $\varsigma=1-a$, $W_0=(2+10\varsigma)I_1+(2-10\varsigma)I_2$ and $Z=-\nu-\E(Fp_1p_2)$. Then
\begin{equation}\label{eq:d2Tabary}
 \Phi(f)=-2Q_L=W_0-8\varsigma Z-2TT-4a\Lambda_1-4b\Lambda_2-2K_F .
\end{equation}
\end{proposition}
\begin{proof}
Cartesian integration by parts gives $\E(p_2\partial_1F)=\E(p_1\partial_2F)=Z$. Since $Lp=(\partial_1F-p_2,\partial_2F+p_1)$, angular integration by parts on $p_1p_2$ gives $2Z+I_1-I_2=-\E(Fp_1p_2)$. Hence $\E(Fp_1p_2)=I_1-I_2-2\nu$ and $Z=\nu-I_1+I_2$. Substitution in \eqref{eq:phi} and expansion of $S_M$ give \eqref{eq:d2Tabary}.
\end{proof}

\subsection{An exact decomposition into squares}
For the square decomposition, the direction tests and the small-temperature theorem, assume also $\E v=0$. Let $R_0(v_1,v_2)=(-v_2,v_1)$, $P=\E(Fp)$, $C_s=\E(p_1p_2)$ and
\[
 \boldsymbol Y=\nabla F+\tfrac12R_0p+\tfrac{J+1}{2}R_0\Theta^{-1}v+P,\qquad
 \kappa_F=\E[(LF+J)^2],\quad \sigma_F^2=\E[(F^2-J)^2] .
\]
Write $(M_1,M_2)=(a,b)$, $A_i=M_i^2\E[(h_i+I_i)^2]$, $r_i=2M_i(\E Y_i^2+P_i^2)$ and $e_H(z)=\E[(H_{12}+C_s-zF)^2]$.

\begin{lemma}\label{lem:d2decomp}
For every $z\in\R$,
\begin{equation}\label{eq:d2decomp}
 Q_L=B_{\rm sc}+\kappa_F+A_1+A_2+r_1+r_2+2ab\bigl(2z\nu-z^2J+C_s^2+e_H(z)\bigr),
\end{equation}
where $B_{\rm sc}=a^2I_1^2+b^2I_2^2+J^2+(J+1)^2+(11a/2-7)I_1+(11b/2-7)I_2$.
\end{lemma}
\begin{proof}
Weighted integration by parts gives $\E F=0$, $\E LF=-J$ and $\E H=-\E(pp^{\mathsf T})$. Thus $K_F=J^2+\kappa_F$, and expansion of the centred Hessian squares gives $TT=a^2I_1^2+b^2I_2^2+A_1+A_2+2ab(2z\nu-z^2J+C_s^2+e_H(z))$. For $U=\nabla F+R_0p/2$ and $V=R_0\Theta^{-1}v$, the midpoint identity for $\nabla F,Lp=\nabla F+R_0p$ gives $S_M=2\E|U|_M^2+(bI_1+aI_2)/2$. We have $\E|V|_M^2=2$, $\E V=0$ and $\E\boldsymbol Y=0$, since $\E(\nabla F)=-P$ and $\E p=\E v=0$. Also
\[
 \E(v_2Y_1)=\E(v_2\partial_1F)-J/2
 =\E(v_1\partial_2F)+J/2=\E(v_1Y_2)
\]
by $\E LF=-J$. Hence $\E(\boldsymbol Y^{\mathsf T}MV)=0$. Squaring $U=\boldsymbol Y-(J+1)V/2-P$ now gives $2\E|U|_M^2=r_1+r_2+(J+1)^2$. Adding the three identities proves the formula.
\end{proof}

\subsection{Pairings fixed by integration by parts}
Put $t_1=b$, $t_2=a$, $D_i=p_i+v_i/t_i$ and $\delta_i=\E D_i^2=I_i-1/t_i\ge0$. Here $\E(v_iD_i)=0$ since $\E(v_ip_i)=-1$. Let $m=\E(v_2\partial_1F)-J/2$ and
\[
 s_1=\nu-I_1+I_2/2+(J+1)/(2a),\qquad s_2=\nu-I_1/2+I_2-(J+1)/(2b).
\]

\begin{lemma}\label{lem:d2pairings}
For $i=1,2$, $k=3-i$ and $\alpha,\beta,\sigma\in\R$,
\begin{align*}
 \E Y_i^2&\ge 2\alpha m+2\beta(s_i+m/t_k)-\alpha^2t_k-\beta^2\delta_k,\\
 a(2\sigma\nu-\sigma^2\delta_2)-aI_1J
 &\le r_1/2+(a\delta_2+J^2)/4+a s_1+(J+1)m+\sqrt{A_1\sigma_F^2},\\
 \kappa_F+A_1-\tfrac25\sqrt{A_1\sigma_F^2}&\ge-4J^2/45 .
\end{align*}
\end{lemma}
\begin{proof}
Integration by parts gives $\E(v_kY_i)=m$ and $\E(p_kY_i)=s_i$. Expand $\E[(Y_i-\alpha v_k-\beta D_k)^2]\ge0$ for the first bound. For the second, set $G=\partial_1F+Fp_1$. Integrating $\partial_1(fF)D_2$ gives $\E(GD_2)=-\nu$, and integrating $\partial_1(fF^2p_1)$ gives $\E G^2=\Lambda_1-\E(F^2h_1)$. The radial cutoff used in \eqref{eq:cubictest} applies to both tests since $F,p_1\in L^4(f)$ and $D_2\in L^2(f)$. Hence $2\sigma\nu-\sigma^2\delta_2\le\E G^2\le\Lambda_1+I_1J+\sqrt{A_1\sigma_F^2}/a$. Expanding $\partial_1F=Y_1+p_2/2+(J+1)v_2/(2a)-P_1$ gives $a\Lambda_1=r_1/2+(a\delta_2+J^2)/4+a s_1+(J+1)m$. If $\delta_k=0$ then $D_k=0$ in $L^2(f)$, so $s_i+m/t_k=0$, and $\delta_2=0$ also gives $\nu=0$. If $J=0$ then $F=0$ in $L^2(f)$ and again $\nu=0$.

The angular cubic test gives $\E[(LF+J)(F^2-J)]=(2J^2-\sigma_F^2)/3$. Cauchy--Schwarz yields $(\sigma_F^2-2J^2)^2\le9\kappa_F\sigma_F^2$. For $\sigma_F^2>0$,
\[
 \kappa_F-\sigma_F^2/25+4J^2/45\ge
 \frac{(\sigma_F^2-2J^2)^2}{9\sigma_F^2}-\frac{\sigma_F^2}{25}+\frac{4J^2}{45}
 =\frac{(4\sigma_F^2-10J^2)^2}{225\sigma_F^2}\ge0 .
\]
For $\sigma_F^2=0$, the same lower bound follows from $\kappa_F\ge0$. Adding $(\sqrt{A_1}-\sqrt{\sigma_F^2}/5)^2$ proves the last bound.
\end{proof}

\subsection{Small cold temperature}
The cold estimate uses the free parameter $z$ in \eqref{eq:d2decomp} to cancel the coefficient of $\nu$ after the direction tests, leaving a polynomial lower bound for $a^2Q_L$.
\begin{theorem}\label{thm:d2small}
If $0<a\le3/5$, then $Q_L>1/(16a^2)$ and $\Phi(f)<-1/(8a^2)$.
\end{theorem}
\begin{proof}
Discard $A_2+2ab(C_s^2+e_H(z))\ge0$ in \eqref{eq:d2decomp}. Apply the second bound of \cref{lem:d2pairings} with multiplier $2/5$ and the third to absorb the square root. This leaves $8a/5$ and $2b$ as the coefficients of $\E Y_1^2$ and $\E Y_2^2$. Set
\[
 r=19/5-3a,\quad t=14/5-3a,\quad h=2a-1,\quad \lambda=r+2a/5,\quad c_w=\lambda-bh-t .
\]
Use $z=c_w/(4ab)$, $\sigma=5t/(4a)$ and choose $(\alpha,\beta)=(m/a,-5r/(16a))$ for $Y_1$, $(m/b,h/4)$ for $Y_2$. The coefficients of $s_1,s_2,\nu$ are $-\lambda,bh,4abz+4a\sigma/5=\lambda-bh$. Since both $s_i$ contain $\nu$ with coefficient one, only $s_i-\nu$ remain. Dropping the nonnegative $P_i^2$ terms and writing $U_2=a/10+5(4t^2+r^2)/(32a)$ gives
\[
\begin{aligned}
 Q_L\ge{}&B_{\rm sc}-\tfrac{17}{90}J^2-\tfrac25 aI_1J
        -U_2\delta_2-\tfrac18 bh^2\delta_1-\tfrac{c_w^2}{8ab}J\\
 &-\lambda(s_1-\nu)+bh(s_2-\nu)+\tfrac{18}{5}m^2
       -\bigl[r/a-h+\tfrac25(J+1)\bigr]m .
\end{aligned}
\]
Put $X=a^2I_1$, $Y=aI_2\ge1$ and $j=aJ$. Multiplying by $a^2$ and completing the square in $am$ leaves the nonnegative term $\tfrac{18}{5}[am-\tfrac5{36}(\lambda-ah+2j/5)]^2$ and the correction $-\tfrac5{72}(\lambda-ah+2j/5)^2$.

After inserting $s_1,s_2,\delta_1,\delta_2$, the terms in $X,j$ are $X^2-2Xj/5+9j^2/5+\ell_X X+\ell_j j$. Completing first the $X$ square and then the $j$ square gives
\[
\begin{aligned}
 a^2Q_L\ge{}&\left(X-\tfrac j5+\tfrac{\ell_X}2\right)^2+
 \tfrac{44}{25}\left[j+\tfrac{25}{88}(\ell_j+\ell_X/5)\right]^2\\
 &+(Y-1)\bigl[b^2(Y-1)+d_Y\bigr]+M_{\rm bd}(a),\\
 M_{\rm bd}(a)={}&N_C/360-\ell_X^2/4-\tfrac{25}{176}(\ell_j+\ell_X/5)^2 .
\end{aligned}
\]
Here $d_Y=[135+455a+a(3-5a)(64a+177)]/160>0$, $\ell_X=N_X/40$, $\ell_j=N_j/(1800b)$, and
\[
\begin{aligned}
 N_X&=20a^3-20a^2+61a-98,\\
 N_j&=-900a^4+5740a^3-17661a^2+24010a-9625,\\
 N_C&=80a^4-1420a^3+2017a^2-1784a+1079 .
\end{aligned}
\]
Write $N_j+9bN_X=b\mathcal Q-729$, with $\mathcal Q=1080a^3-4120a^2+10330a-5330$. Differentiating gives
\[
 22809600\,M_{\rm bd}''=Z_0+2916\mathcal Q'/b^2+1458\mathcal Q''/b
                         +2916\mathcal Q/b^3-6\cdot729^2/b^4,
\]
where $Z_0=-77760000a^4+235008000a^3-532081920a^2+175967040a-100127008$. On $[0,3/5]$, $\mathcal Q''=6480a-8240<0$, $0<\mathcal Q'(3/5)\le\mathcal Q'\le10330$ and $\mathcal Q\le\mathcal Q(3/5)=-9548/25$. Using $a^3\le3a^2/5$ gives $Z_0\le10^6(-390a^2+176a-100)<-80\cdot10^6$, since $390(390a^2-176a+20)=(390a-88)^2+56$. Also $2916\mathcal Q'/b^2\le753057000/49<16\cdot10^6$ since $b\ge7/5$. Thus $22809600\,M_{\rm bd}''<-64\cdot10^6$.
At the endpoints, substitution gives $M_{\rm bd}(0)=1367563/18247680>1/16$ and $M_{\rm bd}(3/5)=74933805223/139708800000>1/2$. Concavity gives $a^2Q_L\ge M_{\rm bd}(a)>1/16$.
\end{proof}

\subsection{Moderate cold temperature}
\begin{theorem}\label{thm:d2moderate}
Assume the hypotheses of \cref{def:staticclass} in dimension two with the zero-mean condition omitted. If $3/5\le a\le1$, then $\Phi(f)\le-J^2/2$. On $3/5\le a\le2/3$ one has $\Phi(f)<-J^2/2-1/(200a^2)$.
\end{theorem}
\begin{proof}
Let $Q_i=\E p_i^4$, $Q_F=\E F^4$, $N=\E H_{12}^2$, $A_R=LF+(F^2+J)/2$ and $B_i=h_i+p_i^2$. The fourth-moment bound in \eqref{eq:cubictest} puts $A_R,B_1,B_2$ in $L^2(f)$. Their pairings are
\[
\begin{aligned}
 \E A_R^2&=K_F-Q_F/12-J^2/4, & \E B_i^2&=\E h_i^2+Q_i/3,\\
 \E(A_RB_1)&=\Lambda_1-2Z, & \E(A_RB_2)&=\Lambda_2+2Z,\\
 \E(B_1B_2)&=2N-\E[(h_1+I_1)B_2].
\end{aligned}
\]
For the first row expand the squares and use \eqref{eq:cubictest}, $\E LF=-J$ and $\E h_i=-I_i$.
For the second row set $F_i=\partial_i F$, $j=3-i$, $e_1=1$ and $e_2=-1$. Then $[\partial_i,L]=e_i\partial_j$. For a radial cutoff $\chi_R$ with $\|\nabla\chi_R\|_\infty\le C/R$, integrate $\partial_i(f\chi_Rp_iA_R)-L(f\chi_Rp_iF_i)$. The third derivatives cancel, and $L\chi_R=0$, giving
\[
 \E\!\left[\chi_R\{B_iA_R-F_i^2+e_i(p_jF_i+p_iF_j)\}\right]
       =-\E[(\partial_i\chi_R)p_iA_R].
\]
The right side tends to zero by Cauchy--Schwarz. Taking $\chi_R\to1$ and using $\E(p_j F_i)=\E(p_i F_j)=Z$ gives the two signs in the second row.

For the last row write $k=H_{12}$. Integrating $\partial_1(f\chi_R p_1h_2)-\partial_2(f\chi_R p_1k)$ cancels the third derivatives because $\partial_1h_2=\partial_2k$. Both boundary coefficients $p_1h_2,p_1k$ lie in $L^1(f)$. Passing to the limit yields $\E(h_1h_2-k^2)=\E(p_1p_2k-p_1^2h_2)$. Integrating $\partial_i(f\chi_R p_i p_j^2)$ in the two coordinate orders gives $\E(h_i p_j^2+2p_i p_jk+p_i^2p_j^2)=0$ for $i\ne j$. Combining these identities yields $\E[(2h_1+p_1^2)B_2]=2N$. Since $\E B_2=0$, this gives the last row.
\emph{Case $0\le\varsigma\le1/3$.} Expanding the square with $\E(A_RB_2)=\Lambda_2+2Z$ gives
\[
 -2K_F-8\varsigma Z+J^2/2=-2\E[(A_R+\varsigma B_2)^2]
 +4\varsigma\Lambda_2+2\varsigma^2\E B_2^2-Q_F/6 .
\]
In \eqref{eq:d2Tabary} the remaining coefficient of $\Lambda_2$ is $-4(b-\varsigma)=-4$. Since $\E B_2^2=\E h_2^2+Q_2/3\le4\E h_2^2$ and $b^2-4\varsigma^2=a(1+3\varsigma)>0$, using $\E h_i^2\ge I_i^2$ gives
\[
 \Phi+J^2/2\le W_0-2a^2 I_1^2-2a(1+3\varsigma) I_2^2.
\]
Write $P_1(x)=(2+10\varsigma)x-2a^2x^2$. If $3\varsigma^2+10\varsigma-1\le0$, then $bP_1'(1/b)=2(3\varsigma^2+10\varsigma-1)\le0$ and $ab^2[16\varsigma/a-P_1(1/b)]=8\varsigma^2(5+3\varsigma)\ge0$. Otherwise $\varsigma>1/19$, since the increasing quadratic is negative at $1/19$. Its sign and $\varsigma\le1/3$ give $(1+5\varsigma)^2-32a\varsigma=(3\varsigma-1)(19\varsigma-1)\le0$. The unrestricted maximum of $P_1$ is therefore at most $16\varsigma/a$. The second score quadratic has derivative $-2-22\varsigma<0$ at $1/a$ and decreases thereafter. Its value at $1/a$ is $-16\varsigma/a$, so the sum is nonpositive.

\emph{Case $1/3\le\varsigma\le2/5$.} Put $\xi=a/30$, $r=\varsigma-\xi$, $W=299(31\varsigma-1)^2/268800$ and $c=b^2-W$. The same moment $\E[(h_1+I_1)B_2]$ enters the two squares
\[
 \E[(A_R-\xi B_1+rB_2)^2],\qquad \E[(a(h_1+I_1)-15rB_2/448)^2] .
\]
Give these squares weights $2,448/225$. In the first weighted square the $A_R B_i$ pairings contribute $8(\xi+r)Z=8\varsigma Z$, exactly the term in the upper bound minus $\Phi+J^2/2$. Their coefficient of $\E[(h_1+I_1)B_2]$ is $4\xi r-2ar/15=0$, and their combined coefficient of $\E B_2^2$ is $(2+1/448)r^2=2W$. For the two remaining squares, \eqref{eq:cubictest} gives
\[
\begin{aligned}
 \E(3h_1+p_1^2+2I_1)^2&=9\E h_1^2-Q_1-4I_1^2,\\
 \E(h_2+p_2^2/23+22I_2/23)^2&=\E h_2^2-43Q_2/1587-484I_2^2/529.
\end{aligned}
\]
Weighting these by $2\xi^2/3,2c$ gives
\[
 \Phi+J^2/2\le W_0-\tfrac{1346}{675}a^2 I_1^2-\tfrac{968}{529}c I_2^2 .
\]
The difference between the right and left sides is the sum of the four weighted square expectations above and
\[
 4(a+\xi)\Lambda_1+4(1+\xi)\Lambda_2+4(ab+2\xi r)N
 +Q_F/6+\tfrac{2}{1587}(43b^2-572W)Q_2 .
\]
Here $r=(31\varsigma-1)/30\ge0$. The ratio $(31\varsigma-1)/(1+\varsigma)$ has derivative $32/(1+\varsigma)^2$ and reaches $57/7$ at $\varsigma=2/5$. Thus $W/b^2\le323817/4390400<43/572<1$. This proves $c>0$ and positivity of the last coefficient.

The $I_2$ quadratic decreases for $I_2\ge1/a$ since $2-10\varsigma<0$ and $c>0$. Take its value at $1/a$ and the unrestricted maximum in $I_1$. After multiplication by $a^2$ this gives $a^2(\Phi+J^2/2)\le\Upsilon(\varsigma)$, where
\[
 \Upsilon(\varsigma)=\tfrac{675}{1346}(1+5\varsigma)^2+(2-10\varsigma)(1-\varsigma)-\tfrac{968}{529}c .
\]
Its quadratic coefficient is $16875/1346+10-968/529+(968\cdot299\cdot961)/(529\cdot268800)>0$. Direct substitution gives $\Upsilon'(1/3)=2779812077/640830600>0$ and $\Upsilon(2/5)=-861265439/99684760000<-1/200$. Thus $\Upsilon$ increases throughout the interval and the right endpoint gives the claimed strict bound. At the common endpoint $a=3/5$ the cold lower bound satisfies $M_{\rm bd}(3/5)>1/2$, while $1-a=2/5$ gives $\Upsilon(2/5)<-1/200$ for the moderate-temperature estimate.
\end{proof}

\begin{theorem}\label{thm:d2}
For every $f\in\cA_2$ one has $\Phi(f)\le0$.
\end{theorem}
\begin{proof}
Interchanging the coordinates exchanges $a,b$ and sends $L$ to $-L$. Thus $\varsigma$ and $Z=\nu-I_1+I_2$ both change sign, while the sums of squares in \eqref{eq:d2Tabary} are preserved. We may therefore take $a\le b$, so $0<a\le1$. The conclusion follows from \cref{thm:d2small,thm:d2moderate}.
\end{proof}

\begin{proof}[Proof of \cref{thm:main}]
Together \cref{thm:d2,thm:d3,thm:dge4} give the static inequality on $\cA_d$ in every dimension $d\ge2$. Applying \cref{thm:reduction} proves the theorem.
\end{proof}

\section{Comparison with the Maxwellian Cauchy problem}\label{sec:discussion}

Tabary~\cite[Theorem~1.1 and the discussion after Remark~1.2]{Tabary} uses the Maxwellian Cauchy solution of Villani~\cite[Proposition~4]{VillaniMaxwell}. To apply our result to that solution, we verify the time trace and the weak equation in \cref{def:weak} with the normalisation of \eqref{eq:landau}.

Let $m>0$, let $c$ be the initial mean velocity, and let $X$ be the continuous Landau process of Delarue, Menozzi and Nualart started with density $f_{\rm in}(c+\cdot)/m$. This probability law has a Lebesgue density and finite second moment, so the non-Dirac hypothesis of their representation holds. Their It\^o formula~\cite[(2.2)--(2.3)]{DMN} has two-particle test kernel $(d-1)K_\varphi/2$ in our notation. Consequently, for $\theta(t)=2mt/(d-1)$, the measures $\mu_t=m\,\mathrm{Law}(c+X_{\theta(t)})$ satisfy
\[
 \langle\mu_t,\varphi\rangle-\langle f_{\rm in},\varphi\rangle
 =\int_0^t\iint K_\varphi(v,w)\,\dd\mu_s(v)\,\dd\mu_s(w)\,\dd s
 \qquad(\varphi\in C_c^\infty).
\]
Indeed the coefficient after the time change is $m\theta'(t)(d-1)/2=m^2$. The representation gives the classical Cauchy density for positive times and conserves mean velocity and total second moment~\cite[Sections~2.1 and~3.3]{DMN}. Taking its classical representative there and a measurable representative of $f_{\rm in}$ at zero gives a jointly measurable density. Continuity of $X$ and dominated convergence give continuity of all $C_c$ pairings up to zero. The bound $|K_\varphi(v,w)|\le C_\varphi(1+|v|^2+|w|^2)$ gives absolute integrability from mass and energy conservation. Integration against a time test now gives \cref{def:weak}. By \cref{prop:evolution}, the Cauchy density agrees with the weak solution almost everywhere in velocity at each time. Its positive-time collision integral is therefore covered by \cref{thm:main}.

Tabary also proves an eventual exponential decay estimate for the entropy production~\cite[Theorem~1.3]{Tabary}. Whether such an estimate holds from time zero remains open.

\preprintdeclarations
\noindent\textbf{Data and code availability.} This theoretical work uses no empirical datasets. A Lean~4 formalisation accompanies the mathematical development.

\smallskip\noindent\textbf{Use of AI tools.} Large language model agents assisted with the mathematical derivations, the Lean formalisation, and the drafting and cross-checking of the manuscript. The authors take responsibility for the content.

\preprintbibliography{KT001}

\end{document}